 \documentclass{ajour}
 \usepackage[all]{xy}
 \usepackage{amsmath,amssymb}
\begin{document}

%


\authorrunninghead{R.~L.~Fernandes}
\titlerunninghead{Lie Algebroids, Holonomy and Characteristic Classes}





\newcommand{\Rr}{\mathbb R}
\newcommand{\Zz}{\mathbb Z}
\newcommand{\norm}[1]{\left\Vert#1\right\Vert}
\newcommand{\abs}[1]{\left\vert#1\right\vert}
\newcommand{\set}[1]{\left\{#1\right\}}
\newcommand{\seq}[1]{\left\langle#1\right\rangle}
\newcommand{\eps}{\varepsilon}
\newcommand{\al}{\alpha}
\newcommand{\e}{\mathbf{e}}
\newcommand{\s}{\mathbf{s}}
\renewcommand{\t}{\mathbf{t}}
\newcommand{\To}{\longrightarrow}
\newcommand{\BX}{\mathbf{B}(X)}
\newcommand{\A}{\mathcal{A}}
\newcommand{\D}{\mathcal{D}}
\newcommand{\F}{\mathcal{F}}
\newcommand{\G}{\mathcal{G}}
\newcommand{\Ho}{\mathcal{H}}
\newcommand{\M}{\mathcal{M}}
\newcommand{\N}{\mathcal{N}}
\newcommand{\X}{\mathcal{X}}
\newcommand{\Lie}{\mathcal{L}}
\renewcommand{\L}{{\text{Lin}}}
\newcommand{\Comp}{\mathcal{K}}
\newcommand{\Basis}{\mathcal{B}}
\newcommand{\Agerm}{\mathfrak{Aut}}
\newcommand{\Ogerm}{\mathfrak{Out}}
\newcommand{\Dgerm}{\mathfrak{Diff}}
\renewcommand{\gg}{\mathfrak{g}}
\newcommand{\hh}{\mathfrak{h}}
\newcommand{\dd}{\mathfrak{d}}
\newcommand{\kk}{\mathfrak{k}}
\newcommand{\mm}{\mathfrak{m}}
\newcommand{\gl}{\mathfrak{gl}}
\newcommand{\Exp}{\text{\rm Exp}\,}
\newcommand{\Ker}{\text{\rm Ker}\,}
\newcommand{\Ad}{\text{\rm Ad}\,}
\newcommand{\ad}{\text{\rm ad}\,}
\newcommand{\Aut}{\text{\rm Aut}\,}
\newcommand{\Out}{\text{\rm Out}\,}
\newcommand{\Inn}{\text{\rm Inn}\,}
\newcommand{\rank}{\text{\rm rank}\,}
\newcommand{\tr}{\text{\rm tr}\,}
\renewcommand{\div}{\text{\rm div}\,}
\title{Lie Algebroids, Holonomy and Characteristic Classes}
\author{Rui Loja Fernandes\thanks{Supported in part by FCT
through program POCTI and grant POCTI/1999/MAT/33081.}}
\affil{Dep.~de Matem\'{a}tica, Instituto Superior T\'{e}cnico,
1049-001 Lisboa, Portugal}
\email{rfern@math.ist.utl.pt}

\abstract{
We extend the notion of connection in order to study
singular geometric structures, namely, we consider a notion of
connection on a Lie algebroid which is a natural extension of the
usual concept of a covariant connection. It allows us to define
holonomy of the orbit foliation of a Lie algebroid and prove a
Stability Theorem. We also introduce secondary or exotic
characteristic classes, thus providing invariants which generalize 
the modular class of a Lie algebroid.}
\begin{article}

\zerosection{Introduction and Basic Definitions}

The theory of connections is a classical topic in differential geometry.
They provide an extremely important tool to study geometric structures on
manifolds and, as such, they have been applied with great success in many
different settings.

However, the use of connections has been very limited whenever
singular be\-ha\-vior is present. The reason is that if some geometric
structure admits a compatible connection then parallel transport
will preserve any algebraic invariant of the structure, and that
prevents the presence of singular behavior. For example, a Poisson
manifold admitting a connection compatible with the Poisson tensor
must have constant rank and hence is a regular Poisson manifold 
(see \cite{Vaisman:book:1}).
In this work we explain how one can  extend the notion of
connection in order to include geometric structures that may
exhibit singular behavior.

One of the basic ideas underlying our construction of connections is
that one should replace the tangent bundle of the manifold $M$ by a
new bundle which reflects more faithfully the (possible singular)
geometric structure on $M$. In this paper we take the point of view
that every such geometric structure has an underlying \emph{Lie
algebroid} structure, which plays the role of the tangent bundle. Many
common geometric structures, some of which we recall below, admit such
a description. We believe that Lie algebroids provide the
appropriate setting to develop a complete geometric theory of
connections as well as other key concepts of differential geometry for
singular geometric structures.

To every Lie algebroid there is associated a foliation of $M$, which in
general will be singular. Conversely, it is not known if every singular foliation 
is associated with a Lie algebroid. There is some evidence that this
is actually the case, and in many ways this work is inspired
by the theory of regular foliations. Several new results to be
presented here are extensions to singular foliations associated with
Lie algebroids of well-known results in foliation theory.

Since the notion of a Lie algebroid is still not part of mainstream
differential geometry, we start by recalling its definition (for an
introduction to the theory see the recent monograph
\cite{Cannas:book}, and also the survey article
\cite{Mackenzie:survey:1}):

\begin{definition}
A \textsc{Lie algebroid} $A$ over a smooth manifold $M$ is a vector
bundle $\pi:A\to M$ together with a Lie algebra structure $[~,~]$
on the space of sections $\Gamma(A)$ and a bundle map $\#:A\to TM$,
called the \emph{anchor}, such that:
\begin{enumerate}
\item[i)] the induced map $\#:\Gamma(A)\to\X^1(M)$ is a Lie algebra
homomorphism(\footnote{We denote by $\Omega^r(M)$ and $\X^r(M)$,
respectively, the spaces of differential $r$-forms and
$r$-multivector fields on a manifold $M$.If $E$ is a bundle over
$M$, $\Gamma(E)$ will denote the space of global sections.});
\item[(ii)] for any sections $\al,\beta\in \Gamma(A)$ and smooth function
$f\in C^\infty(M)$ we have the Leibniz identity:
\begin{equation}
\label{eq:bracket:function}
[\al,f\beta]=f[\al,\beta]+\#\al(f)\beta.
\end{equation}
\end{enumerate}
\end{definition}

The image of $\#$ defines a smooth generalized distribution in $M$,
in the sense of Sussman \cite{Sussmann:article:1}, which is
integrable (this follows, for example, from the Local Splitting Theorem to
be presented below). The integrable leaves are called
\emph{orbits} of $A$ and they form the \emph{orbit foliation} of
the Lie algebroid. We call $A$ a \emph{regular Lie algebroid} if
the rank of $\#$ is locally constant, so the orbit foliation is not
singular. We call $A$ a \emph{transitive Lie algebroid} if $\#$ is
surjective, so the leaves are the connected components of $M$.

The definition of a morphism of Lie algebroids, not necessarily
over the same base manifold, is rather subtle and there are actually two
distinct notions (see e.g.~\cite{Mackenzie:survey:1}). We will be dealing 
mostly with isomorphisms, and for these the two definitions coincide. To 
introduce the definition we will be using, first observe that any 
bundle map $\phi:A_2^*\to A_1^*$ induces a map
$\Phi:\Gamma(A_1)\to \Gamma(A_2)$ which assigns to each section
$\al\in \Gamma(A_1)$ the section $\Phi(\al)\in \Gamma(A_2)$ given
by
\[ \Phi(\al)(y)\equiv\phi^*\al(\phi_0(y)), \qquad \forall y\in M_2,\]
where we have denoted by $\phi_0:M_2\to M_1$ the map induced by
$\phi$ on the base manifolds and  by $\phi^*:(A_1)_{\phi_0(y)}\to
(A_2)_{y}$ the fiberwise transpose of $\phi$.

\begin{definition}
Let $A_1\to M_1$ and $A_2\to M_2$ be Lie algebroids.
A \textsc{morphism of Lie algebroids} from $A_1$ to
$A_2$ is a bundle map $\phi:A_2^*\to A_1^*$ such
that:
\begin{enumerate}
\item[(i)] the induced map $\Phi:\Gamma(A_1)\to \Gamma(A_2)$ preserves
brackets
\begin{equation}
\label{eq:morphism:brackets}
[\Phi(\al),\Phi(\beta)]_2=\Phi([\al,\beta]_1), \qquad \al,\beta\in \Gamma(A_1);
\end{equation}
\item[(ii)] the vector fields $\#_1\Phi(\al)$ and $\#_2\al$ are
$\phi_0$-related:
\begin{equation}
\label{eq:morphism:anchors}
\#_2\al=(\phi_0)_*\#_1\Phi(\al), \quad \al,\beta\in \Gamma(A_1).
\end{equation}
\end{enumerate}
We shall denote such a Lie algebroid homomorphism by $\Phi:A_1\to
A_2$.
\end{definition}

To study global properties of Lie algebroids one needs to consider
connections that are adapted to the orbit foliation.  In this
paper, following the approach in \cite{Fernandes:article:1} for the
special case of Poisson manifolds, we introduce \emph{Lie algebroid
connections} based on the notion of \emph{horizontal lift}. The
basic observation in \cite{Fernandes:article:1}, which also applies in
the present setting, is that one should lift the appropriate geometric
objects rather than tangent vectors as one does in the ordinary
theory of connections.

\begin{definition}
Let $\pi:A\to M$ be a Lie algebroid with anchor $\#:A\to TM$ and
Lie bracket $[~,~]$. An \textsc{$A$-connection} on a fiber bundle
$p:E\to M$ over $M$ is a bundle map $h:p^*A\to TE$, where $p^* A\to
E$ is the pullback bundle of $A$ by $p$, such that the following
diagram commutes:
\[
\xymatrix{
p^*A\ar[r]^{h}\ar[d]_{\widehat{p}}& TE \ar[d]^{p_*} \\
A\ar[r]_{\#} &TM }
\]
\end{definition}

If $(u,\al)\in p^*A$, where $u\in E$ and $\al\in A_x$ with $x=p(u)$,
we call $h(u,\al)$ the \emph{horizontal lift} of $\al$ to the point $u$ in
the fiber over $x$.

Depending on the bundle structure, we may require some additional
conditions on the lift $h$:
\begin{enumerate}
\item[(i)] If $E=P(M,G)$ is a principal bundle with structure
group $G$, then we require $h$ to be $G$-invariant:
\[ h(u\cdot g,\al)=(R_g)_*h(u,\al), \qquad \forall g\in G;\]
\item[(ii)] If $E$ is a vector bundle, then we require $h(u,\cdot)$ to
be linear:
\[ h(u,\al+\beta)=h(u,\al)+h(u,\beta);\]
\end{enumerate}

We now recall some basic classes of Lie algebroids:

\paragraph{Tangent Lie algebroid}
Let $M$ be a manifold. The tangent bundle $TM$ becomes a transitive Lie
algebroid when we take as bracket the usual Lie bracket on vector
fields and as anchor map the identity map. This example is
important for us so we can compare new concepts we shall introduced
for Lie algebroids with the standard ones in ordinary differential
geometry.

If $E\to M$ is a fiber bundle, a $TM$-connection on $E$ is just a
connection in the usual sense. The bundle map $h$ is the horizontal
lift in the usual theory of connections and
$\Ho_u=\set{h(u,\al):\al\in TM}$ is the horizontal distribution (as
in \cite{Kobayashi:book:1}). In this way, we may say that the
theory of $A$-connections is a generalization of the usual theory
of connections. Henceforth, we shall refer to $TM$-connections as
\emph{covariant connections}.

\paragraph{Regular foliations}
Let $A\subset TM$ be an integrable subbundle defining a regular
foliation $\F$ of $M$. A section of $A$ is a vector field in $M$
which is tangent to $\F$. If $X$ and $Y$ are vector fields tangent
to $\F$, their Lie bracket $[X,Y]$ is also a vector field tangent
to $\F$. In this way we have a Lie bracket defined on $\Gamma(A)$
and if we let $\#:A\to TM$ be the inclusion, we obtain a Lie
algebroid. For this Lie algebroid the anchor $\#$ is
\emph{injective} and the orbit foliation is $\F$. Conversely, every
Lie algebroid with anchor map $\#$ injective has a regular
foliation $\F$ and is canonically isomorphic with the Lie algebroid
of $\F$.

If $A$ is the Lie algebroid of a regular foliation $\F$, an
$A$-connection on a fiber bundle $E\to M$ is sometimes called a
\emph{partial connection} along the leaves of $\F$ (see, e.g.,
\cite{Kamber:book:1}). This is because the horizontal lift is only
defined for tangent vectors that are tangent to leaves. Partial
connections were used by Kubarski in \cite{Kubarski:article:1} to
study regular Lie algebroids.

Later we shall see that one can use $A$-connections to define the
$A$-holonomy of the orbit foliation of an \emph{arbitrary} Lie
algebroid $A$. In case $\#$ is injective this holonomy coincides with
the usual holonomy of the theory of regular foliations.

\paragraph{Poisson manifolds}
Let $M$ be a Poisson manifold with Poisson tensor $\Pi\in\X^2(M)$.
Then $\Pi$ determines a bundle map $\#:T^* M\to TM$ as well as a Lie
bracket on the space of differential 1-forms $\Omega^1(M)$, which may
be defined by
\[ [\al,\beta]\equiv\Lie_{\#\al}\beta-\Lie_{\#\beta}\al-\Pi(\al,\beta).\]
It is well known that $(T^* M,[~,~],\#)$ is a Lie algebroid. Some
of the results to be presented in this paper generalize
corresponding results for Poisson manifolds presented in
\cite{Fernandes:article:1}. A $T^*M$-connection on a fiber bundle
$E\to M$ is just a \emph{contravariant connection} on $E$ in the
terminology of \cite{Fernandes:article:1}. This example can be
extended in two distinct directions, namely to Dirac manifolds
(\cite{Courant:article:1}) and to Jacobi manifolds
(\cite{Kerbrat:article:1}).

\paragraph{Transformation Lie algebroids}
Let $\rho:\gg\to \X^1(M)$ be an infinitesimal (right) action of a Lie
algebra on the manifold $M$. The associated \emph{transformation Lie
algebroid} is the trivial bundle $M\times\gg\to M$ with anchor map
$\#:M\times\gg\to TM$ defined by
\[ \#(x,v)\equiv\rho(v)_x,\]
and with Lie bracket
\[ [v,w](x)=[v(x),w(x)]+(\rho(v(x))\cdot w)|_{x}-(\rho(w(x))\cdot v)|_{x},\]
where we identify a section $v$ of $M\times\gg\to M$ with a
$\gg$-valued function $v:M\to\gg$. If $\rho$ can be integrated to a
Lie group action, the orbit foliation of $M$ coincides with the orbits 
of the action. Therefore,
connections on transformations Lie algebroids allows one to study the
singular foliations associated with Lie group actions. For example,
the holonomy (see below) of an orbit is useful in the study of orbit
types.

\paragraph{Bundles of Lie algebras}
Let $A$ be a Lie algebroid with $\#\equiv 0$. For each $x\in M$ one
can define a Lie algebra structure in $A_x$ as follows: if
$\al,\beta\in A_x$ are in the fiber over $x$ choose sections
$\tilde{\al}$ and $\tilde{\beta}$ such that $\tilde{\al}(x)=\al$
and $\tilde{\beta}(x)=\beta$. Then $[\al,\beta]\equiv
[\tilde{\al},\tilde{\beta}](x)$. Using the Leibniz identity one
checks that this definition does not depend on the choice of
sections. Thus we see that such a Lie algebroid is a vector bundle
with varying Lie algebra structure on the fibers. Conversely, any bundle of Lie
algebras determines a Lie algebroid with trivial anchor.

A special case is when $M=\set{*}$, so a Lie algebra can be
consider as a Lie algebroid over a point. It is easy to see that a
linear flat connection is just a representation of a Lie algebra.
\vskip 15 pt

Connections are specially useful to compare geometric structures at
different points of $M$. For non-regular Lie algebroids the orbit
foliation is singular and the dimension of the leaves varies, so
one can only hope to compare spaces at different points of the same
orbit. For that one needs the following fundamental notion of path
which was discovered independently by several authors (e.g., in
\cite{Weinstein:article:4} they are called \emph{admissible paths},
and in \cite{Ginzburg:article:1} they are called \emph{cotangent
paths} in the case $A=T^* M$).

\begin{definition}
\label{defn:A:path}
An \textsc{$A$-path} is a piecewise smooth path $\al:[0,1]\to A$,
such that:
\begin{equation}
\#\al(t)=\frac{d}{dt}\pi(\al(t)), \qquad t\in[0,1].
\end{equation}
The curve $\gamma:[0,1]\to M$ given by $\gamma(t)\equiv\pi(\al(t))$
will be called the \emph{base path} of $\al$.
\end{definition}

Notice that the base path of an $A$-path lies on a fixed leaf of
the Lie algebroid. We shall show that given an $A$-connection one
can define \emph{parallel transport} along any $A$-path. Once the
notion of parallelism is available, one can then proceed to develop
a theory of connections where standard concepts such as
curvature, holonomy, geodesic, etc, make sense. In particular we
show that a linear $A$-connection on a vector bundle $p:E\to M$
gives, in a way entirely analogous to the ordinary case, an
\emph{$A$-derivative} operator
$\nabla:\Gamma(A)\times\Gamma(E)\to\Gamma(E)$ which satisfies:
\begin{enumerate}
\item[(i)] $\nabla_{\al+\beta}\phi=\nabla_\al\phi+\nabla_\beta\phi$;
\item[(ii)] $\nabla_\al(\phi+\psi)=\nabla_\al\phi+\nabla_\al\psi$;
\item[(iii)] $\nabla_{f\al}=f\nabla_\al\phi$;
\item[(iv)] $\nabla_\al(f\phi)=f\nabla_\al\phi+\#\al(f)\phi$;
\end{enumerate}
where $\al,\beta\in\Gamma(A)$, $\phi$, $\psi$ are sections of $E$,
and $f\in C^{\infty}(M)$. Conversely, every such operator is
induced by a linear $A$-connection. Connections from this operational
point of view were first introduced in the
case of Poisson manifolds by Vaisman (\cite{Vaisman:art:1}).
Flat linear $A$-connections on a vector bundle $E$ is an important
special case which has also been studied by several authors usually
under the name of \emph{representations of Lie algebroids}
(\cite{Evens:article:1,Huebschmann:article:1,Xu:article:1,yks:article:1}).

In spite of its formal similarities with ordinary connections,
there are many striking differences in Lie
algebroid connection theory: parallel transport does not depend only
on the base path, the holonomy of a flat $A$-connection may be
non-discrete, etc.

However, just like in ordinary geometry,
$A$-connections are useful to study global properties of Lie
algebroids. Using Lie algebroid connections we show that we have a notion of
\emph{holonomy} of the associated foliation to a Lie algebroid.
For the case of a regular foliation it coincides with the usual
notion of holonomy. We show below that the transversal geometry
to a leaf of a Lie algebroid is described by a germ of an algebroid,
so we have a notion of \emph{transverse Lie algebroid structure}.
The holonomy map is by automorphisms of this transversal
algebroid germ.

In general, holonomy is not homotopy invariant, but factoring out
the inner Lie algebroid automorphisms one obtains a notion of
\emph{reduced holonomy} which is invariant by homotopy, and we can
prove the following analogue of the Reeb Stability Theorem:

\begin{theorem}
Let $L$ be a compact, transversally stable leaf, with finite reduced
holonomy. Then $L$ is stable, i.e., $L$ has arbitrarily small
neighborhoods which are invariant under all inner automorphisms.
Moreover, each leaf near $L$ is a bundle over $L$ with fiber a
finite union of leaves of the transverse Lie algebroid structure.
\end{theorem}

\emph{Linear holonomy} of a Lie algebroid is obtained simply by linearizing
the holonomy homomorphism. In the case of Poisson manifolds,
it was studied by Ginzburg and Golubev in \cite{Ginzburg:article:1}. It
can also be discussed more efficiently from the point of view of
linear Lie algebroid connections and, for each leaf, there is a notion
of \emph{Bott $A$-connection}.

As it was shown in \cite{Evens:article:1}, for a
non-regular Lie algebroid there is a natural vector bundle playing
the role of the normal bundle (over the whole of $M$) to the
singular foliation, and which allows us to introduce the notion of a
\emph{basic connection}: these are linear $A$-connections which
preserve the Lie algebroid structure and restrict in each leaf to
the Bott $A$-connection. Comparing a basic connection to a
riemannian connection, as in the regular theory of foliations, one is
lead to \emph{exotic} or
\emph{secondary Lie algebroid characteristic classes}. These are
$A$-cohomology classes which give information on both the geometry
of the Lie algebroid and the topology of the associated foliation
of $M$. In degree 1, this class actually coincides with the
\emph{modular class} of the Lie algebroid, introduced
by Weinstein in \cite{Weinstein:article:2}.

The remainder of the paper is organized as follows. In Section 1, we
we describe some elementary properties of Lie algebroids and their
differential geometry, including local splitting and transverse
structure which we could not find in the literature. In Section 2, we
sketch the theory of Lie algebroid connections. In Section 3, we
introduce holonomy of a leaf of a Lie algebroid foliation and we prove
the Stability Theorem. In the fourth and final section, we introduce
characteristic classes for Lie algebroids and construct the invariants
we have mentioned above. We also give explicit computations of these
invariants for some classes of Lie algebroids.

Finally, we remark that several authors have considered connections on
Lie algebroids in order to study its global properties (see for
example
\cite{Evens:article:1,Itskov:article:1,Kubarski:article:1,Mackenzie:book:1}).
Until the work of Evans, Lu and Weinstein \cite{Evens:article:1}, all
these authors considered regular, or even
transitive, Lie algebroids. In \cite{Evens:article:1}
connections on non-regular Lie algebroids are used for the first time. 
There, the authors consider zero curvature
$A$-connections on vector bundles, which they call
\emph{representations of Lie algebroids}, to construct the modular
class. These results were extended by Xu
\cite{Xu:article:1} to so called BV-algebras, and Huebschmann in
\cite{Huebschmann:article:1,Huebschmann:article:2} developed a
complete algebraic theory. 

\begin{demo}{Remark} 
In the final stages of preparation of this work, I learn of a preprint
by Marius Crainic \cite{Crainic:preprint:1} where an approach to
secondary characteristic classes for representations of Lie algebroids
is proposed (see also the remark at the end of section
\ref{section:modular:class}).  The discussions I have had with him
after the present paper was submitted, were extremely influential in
shaping my view of the subject. The relationship between the two
approaches is explained in \cite{Crainic:preprint:2}.  Our discussions
eventually led to a solution to the problem of integrating Lie
algebroids to Lie groupoids (see \cite{Crainic:Fernandes}), where we
make use of some of the results presented here. In the preprint
\cite{Ginzburg:article:2}, Ginzburg proposes a $K$-theory for Poisson
manifolds and Lie algebroids which is also related to the present
work.
\end{demo}

\section{The Local Structure of Lie Algebroids}

\subsection{The Dual Lie-Poisson Bracket}
\label{sec:automorphisms}
From now on we fix a Lie algebroid $\pi:A\to M$ over $M$ with anchor
map $\#$ and Lie bracket $[~,~]$.  We let $m$ denote the dimension of
$M$ and we let $r$ denote the rank of $A$. We start by recalling the
construction of a canonical Poisson bracket on the dual bundle $A^*$
(see \cite{Cannas:book}, sect.~16.5).

If one fixes local coordinates $(x^1,\dots,x^m)$ over a
trivializing neighborhood $U$ of $M$ where $A$ admits a basis of
local sections $\set{\al^1,\dots,\al^r}$ over $U$, we have
\emph{structure functions} $b^{is},c^{st}_{u}\in C^\infty(U)$
defined by
\begin{align}
\label{structure:functions:anchor}
\#\al^s&=\sum_{i=1}^m b^{si}\frac{\partial }{\partial x^i}, \qquad
(s=1,\dots,r),\\
\label{structure:functions:bracket}
[\al^{s},\al^{t}]&=\sum_{u=1}^r c^{st}_u \al^u, \qquad (s,t=1,\dots,r).
\end{align}
The defining relations for a Lie algebroid translate into certain
p.d.e's involving the structure functions.

One defines a Poisson structure on $A^*$ as follows. Let
$(\xi^1,\dots,\xi^r)$ denote the linear coordinates on the fibers
of $A^*$ associated with the basis of sections
$\set{\al^1,\dots,\al^r}$. The Poisson bracket $\set{~,~}_A$ on
$C^\infty(A^*)$ is defined by:
\begin{align}
\set{x^i,x^j}&=0,\nonumber\\
\set{x^i,\xi^s}&=-b^{si},\\
\set{\xi^s,\xi^t}&=\sum_u c^{st}_u \xi^u.\nonumber
\end{align}
One checks that this bracket is independent of the choice of local
coordinates and basis. Because this bracket is linear on the
fibers, one also calls it the \emph{dual Poisson-Lie bracket} of
$A$.

Let $\al$ be a section of $A$. Then $\al$ defines in a natural way a
function $f_\al: A^*\to \Rr$ which is linear in the fibers. One has
the following properties of the dual Poisson-Lie bracket.

\begin{proposition}
\label{prop:dual:poisson:bracket}
The assignment $\al\mapsto f_\al$ defines a Lie algebra homomorphism
$(\Gamma(A),[~,~])\to (C^\infty(A^*),\set{~,~}_A)$. Moreover,
if $X_{f_\al}$ denotes the hamiltonian vector field associated
with $f_\al$, then $X_{f_\al}$ is $\pi$-related to $\#\al$:
\[ \pi_*X_{f_\al}=\#\al,\]
where $\pi: A^*\to M$ is the natural projection.
\end{proposition}

\begin{proof}
We use local coordinates. If $\al=\sum_{s}a_s(x)\al^s$ then
$f_\al(x,\xi)=\sum_s a_s(x)\xi^s$ and the associated hamiltonian
vector field is
\begin{equation}
\label{eq:hamiltonian:vector:field}
X_{f_\al}=\sum_{s,i} a_s b^{si}\frac{\partial}{\partial x^i}+
\sum_{t,u}\left(
\sum_{s} a_s c^{st}_u  -
\sum_{i}\frac{\partial a_u}{\partial x^i}b^{t i}
\right)\xi^u \frac{\partial}{\partial \xi^t}.
\end{equation}
This expression shows that $X_{f_\al}$ projects to $\#\al$,

On the other hand, if $\beta=\sum_{t}b_t(x)\al^t$, one computes:
\begin{align*}
\set{f_\al,f_\beta}&=\set{\sum_{s}a_s(x)\xi^s,\sum_{t}b_t(x)\xi^t}\\
                   &=\sum_{t}\#\al(b_t)\xi^t-\sum_{s}\#\beta(a_s)\xi^s
                   +\sum_{s,t,u}a_s b_t c^{st}_u\xi^u=f_{[\al,\beta]},
\end{align*}
so the result follows.
\end{proof}

The dual Lie-Poisson structure of a Lie algebroid codifies all the
information regarding the Lie algebroid structure. In fact, the
category of vector bundles with Poisson brackets linear on the
fibers is equivalent to the category of Lie algebroids. For
example, a morphism $\Phi:A_1\to A_2$ of Lie algebroids, as defined
above, is just a bundle map $\phi:A_2^*\to A_1^*$ which is a
Poisson map.

Let $\al\in\Gamma(A)$ be a section, so we have the associated
hamiltonian vector field $X_{f_\al}$ in $A^*$. For each $t$, the
flow $\phi_t^\al$ defines a Poisson automorphism of $A^*$ (wherever
defined). From (\ref{eq:hamiltonian:vector:field}) we see that
$X_{f_\al}$ is linear along the fibers, so in fact
$\phi_t^\al:A^*\to A^*$ is a bundle map. It follows that each
section determines a 1-parameter family of (local) Lie algebroid morphisms
$\Phi_t^\al:A\to A$. If $A=\gg$ is a Lie algebra, considered as a
Lie algebroid over a point, $\Phi_t^\al: A\to A$ (resp.,
$\phi_t^\al:A^*\to A^*$) is just $Ad(\exp(t\al))$ (resp.,
$Ad^*(\exp(t\al)$), so this construction generalizes the usual
adjoint action, and will be refer to as \emph{integration of
sections}. We remark that we can also integrate time-dependent
sections $\al_t$.

Let us denote by $\Aut(A)$ the group of automorphisms of the Lie
algebroid $A$, and by $\Aut^0(A)$ its connected component of the
identity: given $\Phi\in \Aut^0(A)$ there exists a smooth family
$\Phi_t\in\Aut(A)$, $t\in[0,1]$, such that $\Phi_0=$id,
$\Phi_1=\Phi$. An element $\Phi\in \Aut^0(A)$ is called a
\emph{inner automorphism} if there exists some smooth family of
sections $\al_t\in\Gamma(A)$ which can be integrated to a
1-parameter family of Lie algebroid automorphisms $\Phi_t^{\al_t}$
with $\Phi_1^{\al_t}=\Phi$.

\begin{proposition}
The set $\Inn(A)\subset\Aut(A)$ of inner Lie algebroid
automorphisms is a normal subgroup.
\end{proposition}

We define the group of \emph{outer Lie algebroid automorphisms} to
be the quotient $\Out(A)\equiv \Aut(A)/\Inn(A)$.

\subsection{Local Splitting}
By choosing appropriate coordinates and sections one can simplify
the expressions of the structure functions, and we obtain the
following analogue of the Weinstein Splitting Theorem for Poisson manifolds
(\cite{Weinstein:article:1}, Thm.~2.1).

\begin{theorem}[Local Splitting]
\label{thm:splitting}
Let $x_0\in M$ be a point where $\#_{x_0}$ has rank $q$. There exist
coordinates $(x^i,y^j)$, $(i=1,\dots,q,j=q+1,\dots,m)$, valid in a
neighborhood $U$ of $x_0$, and a basis of sections $\set{\al^1,\dots,\al^r}$,
of $A$ over $U$, such that:
\begin{align}
\label{eq:canonical:anchor:1}
\#\al^i&=\frac{\partial }{\partial x^i}, \qquad (i=1,\dots,q),\\
\label{eq:canonical:anchor:2}
\#\al^s&=\sum_j b^{sj}\frac{\partial}{\partial y^j},\qquad
(s=q+1,\dots,r),
\end{align}
where $b^{sj}\in C^\infty(U)$ are smooth
functions depending only on the $y$'s and vanishing at $x_0$:
$b^{sj}=b^{sj}(y^j)$, $b^{sj}(0)=0$. Moreover,
\begin{equation}
\label{eq:canonical:brackets}
[\al^{s},\al^{t}]=\sum_u c^{s t}_u \al^u,
\end{equation}
where $c^{s t}_u\in C^\infty(U)$ vanish if $u\le q$ and satisfy
\begin{equation}
\label{eq:canonical:structure:functions}
\sum_{u>q}\frac{\partial c^{s t}_u }{\partial x^i}b^{uj}=0.
\end{equation}
\end{theorem}

\begin{proof} If the rank of $\#$ at $x_0$ is $q=0$ we are done, so we
can assume $q\ge 1$. If $q\ge 1$ we proceed, by induction,
straightening out vector fields 
of the form $\#\al$. So let $0\le k<q$ and assume we have constructed
coordinates
\[ (x^{i},\tilde{y}^{j}), \text{ where }i\le k,\ k<j\le m,\]
valid on a domain $U$, and a basis of sections for $A$ over $U$,
\[ \set{\al^i,\tilde{\al}^{s}}, \text{ where }i\le k,\ k<s\le r,\]
such that
\begin{align*}
\#\al^i &=\frac{\partial }{\partial x^i}, \qquad (i\le k),\\
\#\tilde{\al}^s &=\sum_j b^{sj}\frac{\partial}{\partial \tilde{y}^j},\qquad
(k<s\le r),
\end{align*}
where $b^{sj}\in C^\infty(U)$ depend only on the $\tilde{y}$'s. Since
$q>k$, there exists an $s$
such that the vector field $\#\tilde{\al}^s$ does not vanish at $x_0$. By
relabeling, we can assume that $s=k+1$ and we set
$\al^{k+1}=\tilde{\al}^{k+1}$.

By straightening out $\#\al^{k+1}$, we can perform a change of
coordinates
\begin{align*}
x^{k+1}&=x^{k+1}(\tilde{y}^{k+1},\dots,\tilde{y}^m),\\
y^j&=y^j(\tilde{y}^{k+1},\dots,\tilde{y}^m), \qquad j=k+2,\dots,m,
\end{align*}
such that
\begin{align*} 
\#\al^{k+1}&=\frac{\partial }{\partial x^{k+1}}, \\
\#\tilde{\al}^s&=b^{s,k+1}\frac{\partial}{\partial x^{k+1}}+\cdots 
\end{align*}
Replacing $\tilde{\al}^s$ by $\tilde{\al}^s-b^{s,k+1}\al^{k+1}$, we see
that we can assume $b^{s,k+1}=0$. Therefore,
\[ \#\tilde{\al}^s=\sum_{j}b^{sj}\frac{\partial}{\partial y^j},\]
where $b^{sj}=b^{sj}(x^{k+1},y^{k+2},\cdots,y^m)$.

Using $\#[\al^{k+1},\tilde{\al}^s]=[\#\al^{k+1},\#\tilde{\al}^s]$
for $s>k+1$, we see that
\[ [\al^{k+1},\tilde{\al}^s]=\sum_{t>k+1}c^{s,k+1}_t \tilde{\al}^t,\]
where the structure functions are related by
\[ \frac{\partial b^{sj}}{\partial x^{k+1}}=\sum_{u>k+2} c^{s,k+1}_u
b^{uj}.\]
We can think of this equation as a time-dependent linear o.d.e.~for
$b^{sj}$ in the variable $x^{k+1}$. Let us denote by $X(x^{k+1})$ the
fundamental matrix of solutions such that $X(0)=I$, and by
$Y(x^{k+1})$ its inverse. We consider new sections
\[ \al^s=\sum_{t>k+2} Y^s_t(x^{k+1})\tilde{\al}^t.\]
Then we find
\begin{align*}
\#\al^s&=\sum_j\sum_{t>k+2}Y^s_t(x^{k+1})b^{tj}(x^{k+1},y^{k+2},\cdots,y^m)
\frac{\partial}{\partial y^j}\\
       &=\sum_j b^{sj}(0,y^{k+2},\cdots,y^m)\frac{\partial}{\partial y^j}.
\end{align*}
We conclude that there exist coordinates $(x^i,y^j)$ and sections
$\set{\al^s}$, as in the statement of the theorem, such that
(\ref{eq:canonical:anchor:1}) and (\ref{eq:canonical:anchor:2}) hold,
for some smooth functions $b^{sj}\in C^\infty(U)$ depending only on
the $y$'s. Since at $x_0$ the bundle map $\#$ has rank $q$, we must
have $b^{sj}(0)=0$.

Comparing coefficients of $\frac{\partial}{\partial x^i}$ in
$\#[\al^s,\al^t]=[\#\al^s,\#\al^t]$ we check easily that the
structure functions $c^{s t}_u\in C^\infty(U)$ must vanish for
$u\le q$. Using the Jacobi identity, we find for $i\le q$ and
$q<s,t\le r$,
\[\#[\al^i,[\al^s,\al^t]]=[[\#\al^i,\#\al^s],\#\al^t]]+
[\#\al^s,[\#\al^i,\#\al^t]]=0.\]
On the other hand,
\begin{align*}
\#[\al^i,[\al^s,\al^t]]&=\#[\al^i,\sum_{u>q} c^{st}_u \al^u]\\
                       &=[\#\al^i,\sum_{u>q} c^{st}_u\#\al^u]
                       =\sum_j\sum_{u>q} \frac{\partial
                       c^{st}_u}{\partial x^i}
                       b^{uj}\frac{\partial}{\partial x^j},
\end{align*}
so (\ref{eq:canonical:structure:functions}) follows.
\end{proof}

In general, the structure functions that appear in relations
(\ref{eq:canonical:brackets}) will depend both on the $x$'s and
$y$'s variables, subject to
(\ref{eq:canonical:structure:functions}). For special classes of
Lie algebroids one might have extra information that leads to
further simplification of the structure functions. For example, in
the case of a Poisson manifold, one always has the relationship:
\[ c^{ij}_k=\frac{\partial b^{ij}}{\partial x^k}.\]
Then, all structure functions in (\ref{eq:canonical:brackets})
depend only on the $y$'s variables, and one obtains the
Weinstein Splitting Theorem.

Note that Theorem \ref{thm:splitting} \emph{is not} the
Weinstein Splitting Theorem for the Lie-Poisson structure on $A^*$. The
reason is that we are only allowed to make changes of coordinate of
the base manifold $M$ and of sections of $A$. These lead to changes of
coordinates of $A^*$ which are linear in the fiber variables:
\[y^i=y^i(x^1,\dots,x^m), \qquad \eta^s=\sum_t a^s_t(x^1,\dots,x^m)\xi^t.\]
These changes of coordinate are usually not sufficient to obtain the
Weinstein splitting for $A^*$.

As we mentioned above, a simple consequence of the Local Splitting Theorem is that
the generalized distribution Im$\#$ is integrable.

\subsection{Transverse Structure}
\label{sec:trasnverse:structure}
The local splitting of a Lie algebroid can be used to define a
transverse Lie algebroid structure, similar to the case of a Poisson
manifold. We first give a more invariant description, and later
come back to the local coordinate approach.

First we observe that every (embedded) submanifold $N\subset M$ which
is transverse to the orbit foliation
\[ TN+\text{Im}(\#)=TM,\]
has a natural induced Lie algebroid structure $A_N\to N$. We
take for $A_N$ the vector bundle over $N$ with fibers
\[ (A_N)_x=\set{\al\in A_x: \#\al\in T_x N}.\]
Because of the transversality assumption, this is indeed a subbundle
of $A$. The anchor map $\#:A_N\to TN$ is obtained simply by restriction of
$\#$. Also, every section $\al$ in $\Gamma(A_N)$ extends to a section
$\tilde{\al}$ of $A$ defined in an open set containing $N$, and given
two sections $\al,\beta\in\Gamma(A_N)$, we set
\[ [\al,\beta]_{A_N}(x)\equiv[\tilde{\al},\tilde{\beta}](x).\]
One checks that (i) this bracket does not depend on the extensions
considered, and (ii) that it defines a section of $A_N$. It follows
that $(A_N,\#,[~,~]_{A_N})$ is a Lie algebroid over $N$.

The notion of \emph{transverse Lie algebroid structure} is based on
the following result, also inspired in Poisson geometry 
(cf.~\cite{Weinstein:article:1}, Section 2).

\begin{theorem}[Transverse Structure]
Let $L$ be a leaf of the orbit foliation of $A$, and suppose $N_0$
and $N_1$ are submanifolds of $M$ of complementary dimension to $L$
and intersecting $L$ transversally on a single point. Then there
exists an automorphism of $A$ which maps a neighborhood $V_0$ of
$N_0\cap L$ in $N_0$ onto a neighborhood $V_1$ of $N_1\cap L$ in
$N_1$, and which induces an isomorphism of the induced Lie
algebroid structures on the neighborhoods.
\end{theorem}

\begin{proof}
If $x_0=L\cap N_0$ and $x_1=L\cap N_1$, there exists a piece-wise
smooth path made of orbits of vector fields of the form $\#\al$, with
$\al$ a section of $A$. Integrating sections we can map $x_0$ to
$x_1$, so we may assume that these points of intersection are actually
the same.

Around $x_0=x_1$ we choose coordinates $(x^i,y^j)$ and sections
$\set{\al^s}$ as in the Local
Splitting Theorem. We interpolate between $N_0$ and $N_1$ by a
family of manifolds $N_t$ defined by equations of the
form
\[ x^i=X^i(y^1,\dots,y^{m-q},t), \qquad (i=1,\dots,q).\]
Then we look for a time-dependent section $\al_t$ which, by
integration, gives a Lie algebroid automorphism $\Phi_t:A\to A$,
covering a diffeomorphism $\phi_t:M\to M$, which maps a
neighborhood of $x_0=x_1$ in $N_0$ onto a neighborhood of $N_t$.

Let us write $\al_t=\sum_s a_s(x^i,y^j,t)\al^s$. In order for the
$\phi_t$ to track the $N_t$ we must have the equations
\[ a_i=\sum_{j,s} \frac{\partial X^i}{\partial y^j} b^{sj}a_s+
\frac{\partial X^i}{\partial t},\qquad (i=1,\dots,q)\]
satisfied along $N_t$. It is clear than one can choose $a_s$ such that
this equations holds. Integration of $\al_t$ gives a Lie algebroid
automorphism 
\[ \Phi_t^{\al_t}:A\to A\]
which induces a Lie algebroid isomorphism between $A_{N_0}$ and $A_{N_t}$.
\end{proof}

The transversal geometry to a Lie algebroid, around a point, is described 
by a \emph{transversal algebroid germ}, that is to say a germ of a Lie algebroid
for which the anchor vanishes at the base point. For all cross sections 
$N$ to $L$ the induced Lie algebroid structures $A_N$ are locally isomorphic, 
but there is no natural choice for this transverse structure.
In other words,, we have a well defined notion of \emph{transverse Lie algebroid
structure} along a leaf $L$. 

In the local splitting coordinates $(x^i,y^j)$ and sections
$\set{\al^s}$ given by Theorem \ref{thm:splitting},
the transverse Lie algebroid structure $A_N\to N$, has
coordinates $(y^j)$ in the base, sections
$\set{\bar{\al}^s(y)\equiv\al^s(0,y^j):\ q<s\le r}$, anchor map
\[ \#\bar{\al}^s=\sum_{j>q}b^{sj}(y)\frac{\partial}{\partial y^j}, \qquad
(s>q)\]
and Lie algebra structure
\[ [\bar{\al}^s,\bar{\al}^t]=\sum_{u>q}c^{st}_u(0,y)\bar{\al}^u.\]

Now the Local Splitting Theorem can be stated as
follows: Let $x_0$ be any point in a Lie algebroid $A$, and denote by
$L$ the leaf through $x_0$ and by $N$ a cross-section to $L$ at
$x_0$. Then, locally, $A$ is an extension of $A_N$ by $TL$, i.e.,
we have Lie algebroids morphisms
\[ A_N\to A\to TL\]
such that the corresponding bundle maps form a short exact sequence:
\[ 0\to T^*L\to A^*\to A^*_N\to 0.\]

\subsection{Linear Approximation to a Lie Algebroid}

By a \emph{linear Lie algebroid} we shall mean a Lie algebroid $\pi:A\to M$
satisfying the following properties:
\begin{enumerate}
\item[(i)] The base $M=V$ is a vector space so $\pi:A\to M$ is a
trivial bundle;
\item[(ii)] For any trivialization, the bracket of any constant sections
$\al$ and $\beta$ is a constant section $[\al,\beta]$;
\item[(iii)] For any trivialization, the vector field $\#\al$ is linear whenever
$\al$ is a constant section;
\end{enumerate}
Conditions (i) and (ii) mean that $A$ is a transformation Lie
algebroid: $A=\gg\times V$ and $\#:\gg\to \X^1(V)$ is an action of
the Lie algebra $\gg$ on $V$. Condition (iii) means that this
action is linear. So a linear Lie algebroid is just a
representation of a Lie algebra.

Now let $A$ be any Lie algebroid and fix $x_0\in M$. The normal
space $N_{x_0}=T_{x_0}M/\text{Im}\#_{x_0}$ carries a natural linear action of 
the Lie algebra $\mathfrak{g}=Ker\#_{x_{0}}$, and the associated transformation 
Lie algebroid is called the \emph{linearization of $A$ at $x_{0}$}.
This linearization can be seen in two different ways:
\begin{enumerate}
\item[(a)] We take $\gg=Ker \#_{x_0}$ with the Lie algebra
structure induced from $A$, and we define a linear action of $\gg$ on
$N_{x_0}$ by
\[ \rho(z)=\delta_{x_0}(\#\al),\qquad z\in\gg\]
where $\al$ is such that $\al_{x_0}=z$ and $\delta_{x_0}X$ is the
linearization of the vector field $X$ at $x_0$. Since the flow of
$\#\al$ preserves the orbit foliation this actually defines a
linear endomorphism of $N_{x_0}=T_{x_0}M/\text{Im}\#_{x_0}$. If
$\al'$ is another section with $\al'_{x_0}=z$, one checks that
$\delta_{x_0}(\#\al)$ and $\delta_{x_0}(\#\al')$ induce the same
endomorphism of $N_{x_0}$. The associated transformation Lie
algebroid $\gg\times N_{x_0}$ is the linearization of $A$.
\item[(b)] Again we take the trivial vector bundle $\gg\times
N_{x_0}\to N_{x_0}$ and we define the anchor $\tilde{\#}:\gg\times
N_{x_0}\to T N_{x_0}$ to be the intrinsic derivative at $x_0$ of
the bundle map $\#:A\to TM$. Then we define the bracket on
constant sections to be the pointwise bracket, and we extend to any
section by requiring the Leibniz identity to hold.
\end{enumerate}

If we pick splitting coordinates $(x^i,y^j)$ and
basis of sections $\set{\al^s}$ around $x_0$, then
$\set{e^s\equiv\al^s(x_0):s>q}$ give a basis for $\gg=\Ker \#_{x_0}$,
with structure constants $c^{st}_u(x_0)$. The tangent vectors
$v_j=\left.\frac{\partial}{\partial y^j}\right|_{x_0}$ induce a basis for
$N_{x_0}$, determining linear coordinates $(w^j)$ and
relative to these coordinates the anchor map is given by:
\[ \tilde{\#}(e^s)=\sum_{j,k}
\frac{\partial b^{s j}}{\partial y^k}(0)w^k \frac{\partial}{\partial w^j},\]
where we view the $\set{e^s}$ as constant sections of $\gg\times
N_{x_0}$.

One should notice that the transverse Lie algebroid structure is an
equivalence class of isomorphic structures, for which
there is no natural choice of representative, while the linearization at
$x_0$ lives on a well-defined bundle over the normal space
$N_{x_0}$. The problem of linearizing a Lie algebroid is discussed in
\cite{Weinstein:article:3}.

\subsection{Lie Algebroid Cohomology}
The existence of a Lie bracket on the space of sections of a Lie
algebroid leads to a calculus on its sections analogous to the usual
Cartan calculus on differential forms. In this paragraph we give only
the relevant formulas for Lie algebroid cohomology we shall need
later, and refer the reader to the monograph \cite{Mackenzie:book:1}
for further details.

One defines the exterior differential
$d_A:\Gamma(\wedge^\bullet A^*)\To\Gamma(\wedge^{\bullet+1}A^*)$ by:
\begin{multline}
\label{eq:differential}
d_A Q(\al_0,\dots,\al_r)=\frac{1}{r+1}\sum_{k=0}^{r+1}
(-1)^k\#\al_k(Q(\al_0,\dots,\widehat{\al}_k,\dots,\al_r)\\
+\frac{1}{r+1}\sum_{k<l}(-1)^{k+l+1}Q([\al_k,\al_l],\al_0,\dots,\widehat{\al}_k,\dots,\widehat{\al}_l,\dots,\al_r).
\end{multline}
where $\al_0,\dots,\al_r$ are sections of $A$. This differential
satisfies:
\begin{align}
d_A^2(Q)&=0,\\
\label{eq:derivation}
d_A(Q_1\wedge Q_2)&=d_A Q_1 \wedge Q_2+(-1)^{\deg Q_1} Q_1\wedge
d_A Q_2.
\end{align}
The cohomology associated with $d_A$ is called the
\emph{Lie algebroid cohomology} of $A$ (with trivial coefficients) and
is denoted by $H^\bullet(A)$.

Define a homomorphism of exterior algebras
$\#^*:\Omega^\bullet(M)\to \Gamma(\wedge^\bullet A^*)$ by setting:
\[ \#^*\omega(\al_1,\dots,\al_r)=(-1)^r\omega(\#\al_1,\dots,\#\al_r).\]
We compute
\begin{equation}
\label{eq:differentials}
d_A\#^*\omega=\#^*d\omega.
\end{equation}
so there is a ring homomorphism
\[\#^*:H^\bullet_{\text{de Rham}}(M)\to H^\bullet(A).\]

For the examples we have mentioned above these cohomology groups are well
known. Special cases to be used later are:
\begin{enumerate}
\item[(i)] $A=TM$, where $H^\bullet(A)=H^\bullet_{\text{de Rham}}(M)$ is 
the \emph{de Rham cohomology};
\item[(ii)] $A=T^* M$ with $(M,\Pi)$ a Poisson manifold, where one obtains the
\emph{Poisson cohomology} of $M$ denoted $H^\bullet_{\Pi}(M)$;
\item[(iii)] $A=T\F\subset TM$ an integrable subbundle associated with a
regular foliation $\F$, where one gets the
\emph{tangential cohomology} denoted $H^\bullet_{\F}(M)$;
\item[(iv)] $A=\gg\times V$ the Lie algebroid associated with a Lie
algebra representation $\rho:\gg\to \Aut(V)$, where one gets the
\emph{Chevalley-Eilenberg cohomology} $H^\bullet(\gg,\rho)$;
\end{enumerate}
The Lie algebroid cohomology, in general, will not be homotopy invariant
and hence it may be hard to compute (to say the least). This is intimately
related with the singular behaviour of the orbit foliation, one of the
main topics to be discussed here.

Finally, we note, for later reference, that if $\phi:M\to N$ is a smooth
map its \emph{$A$-differential} is the bundle map
$d_A \phi:A\to TN$ defined by:
\begin{equation}
d_A \phi(\al_x)=d_x \phi\cdot\#\al_x, \qquad \al_x\in A_x.
\end{equation}
If $N=\Rr$ this notation is consistent with the A-differential
introduced above, if we think of sections in $\Gamma(\wedge^0 A^*)$ as
functions on $M$.

\section{Lie Algebroid Connections}

\subsection{Connections on Principal Bundles} Let $P(M,G)$ be a smooth
principal bundle over the manifold $M$ with structure group $G$. We
let $p:P\to M$ be the projection, and for each $u\in P$ we denote
by $G_u\subset T_u(P)$ the subspace consisting of vectors tangent
to the fiber through $u$. If we let $\sigma:\gg\to\X^1(P)$ be the
infinitesimal $G$-action on $P$, we have $G_u=\set{\sigma(B)_u|\ B\in\gg}$.

Now let $p^*A$ denote the pullback bundle of $A$ by $p$.
There is a bundle map $\widehat{p}:p^*A
\to A$ which makes the following diagram commutative
\[
\xymatrix{
p^*A\ar[r]^{\widehat{p}}\ar[d]_{\widehat{\pi}}& A \ar[d]^{\pi}
\\ P\ar[r]_{p} &M }
\]
where on the vertical arrows we have the canonical projections.
Recalling that $p^*A=\set{(u,\al)\in P\times A:p(u)=\pi(\al)}$, we
see that we have a natural right $G$-action on $p^*A$ defined by
$(u,\al)\cdot a\equiv(u a,\al)$, if $a\in G$. Our basic definition is
then the following:

\begin{definition}
An \textsc{$A$-connection} in the principal bundle $P(M,G)$ is a
smooth bundle map $h:p^*A\to TP$, such that:
\begin{enumerate}
\item[(CI)] $h$ is horizontal, i.e., the following diagram commutes:
\[
\xymatrix{
p^*A\ar[r]^{h}\ar[d]_{\widehat{p}}& TP \ar[d]^{p_*} \\ A\ar[r]_{\#}
&TM }
\]
\item[(CII)] $h$ is $G$-invariant, i.e., we have
\[ h(u a,\al)=(R_a)_*h(u,\al), \qquad \text{ for all }a\in G;\]
\end{enumerate}
\end{definition}

The subspace of $T_u P$ formed by all horizontal lifts $h(u,\al)\in
T_u P$, where $(u,\al)\in p^*A$, is denoted by $\Ho_u$. The
assignment $u\mapsto \Ho_u$ is a smooth distribution on $P$ called
the \emph{horizontal distribution} of the connection(\footnote{In this paper
``smooth distributions'' are always in the sense of Sussman 
\cite{Sussmann:article:1}, so that for
each point $u_0\in P$ there exists a neighborhood $u_0\in U\subset
P$ and smooth vector fields $X_1,\dots,X_r$ in $U$, such that
$\Ho_u=\text{span}\set{X_1|_u,\dots,X_r|_u}$ for all $u\in U$.}).
Note that, unlike the ordinary case, the rank of the horizontal
distribution will vary, and that this distribution does not define
the connection uniquely.

It follows from (CI) in the definition of an $A$-connection, that
the horizontal spaces $\Ho_u$ project onto the tangent space $T_x
L$ to the orbit leaf $L$ through $x=p(u)$. In general, we have
neither $T_u P=G_u+\Ho_u$ nor $G_u\cap \Ho_u=\set{0}$. As usual, a
vector $X\in T_u P$ will be called \emph{vertical}
(resp.~\emph{horizontal}), if it lies in $G_u$ (resp.~$\Ho_u$).
Since, in general, a tangent vector to $P$ does not split into a
sum of an horizontal and a vertical component, the usual
definitions of lift of curves, connection form, etc., do not make
sense in this context. We will show below how to define these
notions appropriately.

Recall that the \emph{Atiyah sequence} of the principal bundle
$P(M,G)$ is the short exact sequence of vector bundles
\[ \xymatrix{0\ar[r]& \Ad(P)\ar[r]^j &TP/G \ar[r]^{p_*} &TM\ar[r] &0},\]
where $Ad(P)=\frac{P\times \gg}{G}$ is the associated bundle to $P$
obtained from the adjoint representation of $G$ on $\gg$ (in fact,
this is a short exact sequence of Lie algebroids for the obvious
Lie algebroid structures). Now, if $h:p^* A\to TP$ is a connection,
the $G$-invariance implies that $h$ induces a bundle map
$\omega:A\simeq p^*A/G\to TP/G$. The following commutative diagram,
which was suggested to me by Alan Weinstein, is helpful in
understanding the relationship between the different geometric
objects associated with an $A$-connection:
\[
\xymatrix{
0\ar[dr] & & & & & & 0\ar[dl]\\
 & \Ker\#\ar[dr]\ar@{-->}[rrrr]&&&&\Ad(P)\ar[dl]\\
 &&A\ar[dr]_{\#}\ar[rr]^{\omega}&&TP/G\ar[dl]^{p_*}\\
 &&&TM\ar[d]\\
 &&&0}
\]
Note, however, that in this diagram $\Ker\#$ is not a vector
bundle, except if $A$ is a regular Lie algebroid. The dash arrow
will be explained later.

In the usual theory of covariant connections, one has $A=TM$ and $\#$
is the identity map, so a connection can be thought as a splitting of
the Atiyah sequence. This is the approach taken by Mackenzie in
\cite{Mackenzie:book:1} and which lead him to the introduction of a
connection for a \emph{transitive Lie algebroid} as a splitting of the
analogous short exact sequence:
\[ \xymatrix{0\ar[r]& \Ker\#\ar[r] &A \ar[r]^{\#} &TM\ar[r] &0}.\]
This approach was also followed by Kubarski in his theory of
characteristic classes for regular Lie algebroids
(\cite{Kubarski:article:1,Kubarski:article:1}).

\subsection{Connection 1-section and Curvature 2-section}

Let $E\to M$ be any vector bundle. In the theory of Lie algebroids,
elements of $\Gamma^\bullet(A^*,E)\equiv\Gamma(\wedge^\bullet
A^*)\otimes\Gamma(E)$ play the role of ($E$-valued) differential
forms. We shall refer to an element in $\Gamma^r(A^*,E)$ as an
$E$-valued $r$-section, or simply as an \emph{$r$-section} if it is
clear from the context what $E$ is. In case $E$ is also a Lie
algebroid, we have an induced (super) Lie bracket on
$\Gamma^\bullet(A^*,E)$ by setting:
\begin{multline}
\label{eq:super:Lie:bracket}
[P,Q](\al_1,\dots,\al_{p+q})=\\\frac{1}{(p+q)!}\sum_{\sigma}(-1)^\sigma
[P(\al_{\sigma(1)},\dots,\al_{\sigma(p)}),Q(\al_{\sigma(p+1)},\dots,\al_{\sigma(p+q)})],
\end{multline}
where the sum is over all permutations $\sigma$ of $p+q$ elements.

To deal with connections we let $E=TP/G$. As we remarked above, a
connection $h$ defines a bundle map $\omega:A\simeq p^*A/G\to
TP/G$, i.e., an element in $\Gamma^1(A^*,TP/G)$, and we call
$\omega$ the
\emph{connection 1-section}. We define the
\emph{exterior $A$-derivative}
\[D:\Gamma^\bullet(A^*,TP/G)\to \Gamma^{\bullet+1}(A^*,TP/G)\]
by setting:
\begin{multline}
\label{eq:exterior:A:derivative}
D Q(\al_0,\dots,\al_q)=
\frac{1}{q+1}\sum_{k=0}^q
(-1)^k[\omega(\al_k),Q(\al_0,\dots,\widehat{\al}_k,\dots,\al_q)\\
+\frac{1}{q+1}\sum_{k<l}(-1)^{k+l}Q([\al_k,\al_l],\al_0,\dots,\widehat{\al}_k,\dots,\widehat{\al}_l,\dots,\al_q).
\end{multline}

Now we introduce the \emph{curvature 2-section} $\Omega\in
\Gamma^2(A^*,TP/G)$ of the connection by setting:
\begin{equation}
\label{eq:curvature:geometric}
\Omega(\al,\beta)\equiv\frac{1}{2}\left([\omega(\al),\omega(\beta)]-
\omega([\al,\beta])\right),
\end{equation}
for $\al,\beta\in\Gamma(A)$. The curvature 2-section
measures to which extent the horizontal distribution fails to be integrable.

\begin{proposition}
The curvature 2-section satisfies the structure equation:
\begin{equation}
\label{eq:first:structure:equation}
\Omega=D\omega-\frac{1}{2}[\omega,\omega].
\end{equation}
\end{proposition}

\begin{proof}
We compute using (\ref{eq:super:Lie:bracket}) and
(\ref{eq:exterior:A:derivative}):
\begin{align*}
D\omega(\al,\beta)=[\omega(\al),\omega(\beta)]-
\frac{1}{2}\omega([\al,\beta]),\\
[\omega,\omega](\al,\beta)=[\omega(\al),\omega(\beta)]
\end{align*}
so we have
\[ \Omega(\al,\beta)=D\omega(\al,\beta)-\frac{1}{2}[\omega,\omega](\al,\beta)=
\frac{1}{2}\left([\omega(\al),\omega(\beta)]-\omega([\al,\beta])\right),\]
which shows that (\ref{eq:first:structure:equation}) is satisfied.
\end{proof}

The horizontal distribution in general will have non-constant
rank. Still, if we call a \emph{flat $A$-connection} an $A$-connection
for which the horizontal distribution is integrable, we have:

\begin{corollary}
An $A$-connection is flat iff its curvature 2-section vanishes.
\end{corollary}

\begin{proof}
By a result of Hermann \cite{Hermann:article:1}, a generalized
distribution associated with a vector subspace $\D\subset\X^1(M)$
is integrable iff it is involutive and rank invariant. Taking
$\D=\set{h(\al):\al\in \Gamma(A)}$, so that
$\Ho_u=\set{X(u):X\in\D}$, (\ref{eq:curvature:geometric}) shows
that $\D$ is involutive iff the curvature 2-section vanishes.
Hence, all it remains to show is that if the curvature vanishes and
$\gamma(t)$ is an integral curve of $h(\al)$ then
$\dim\Ho_{\gamma(t)}$ is constant, for all small enough $t$.

Fix $\al\in\Gamma(A)$ and let $\phi_t^\al$ be the flow $X_{f_\al}$,
let $\psi_t^\al$ be the flow of $\#\al$ and let $\tilde{\psi}_t^\al$
be the flow of $h(\al)$. We have
$\psi_t^\al=p\circ\tilde{\psi}_t^\al=\pi\circ \phi_t^\al$ (see
prop.~\ref{prop:dual:poisson:bracket}).
If $\beta\in\Gamma(A)$ we claim that
\[ (\tilde{\psi}_t^\al)_* h(\beta)=h(\phi_t^\al\beta), \]
for small enough $t$. In fact, the infinitesimal version of this
relation is
\[ [h(\al),h(\beta)]=h([\al,\beta]), \]
which holds, since we are assuming that the curvature vanishes.

Therefore, the flow $\tilde{\psi}_t^\al$ gives an isomorphism between
$\Ho_{\gamma(0)}$ and $\Ho_{\gamma(t)}$, for small enough $t$, so $\D$
is rank invariant.
\end{proof}

We also have an analogue of the usual Bianchi's identity:
\begin{proposition}
The curvature 2-section $\Omega$ satisfies the Bianchi identity:
\begin{equation}
\label{eq:Bianchi}
D\Omega=0.
\end{equation}
\end{proposition}

\begin{proof}
From expression (\ref{eq:curvature:geometric}) for the curvature
and the definition (\ref{eq:exterior:A:derivative}) of the exterior
$A$-derivative, we compute:
\begin{multline*}
D\Omega(\al,\beta,\gamma)=
\bigodot_{\al,\beta,\gamma}\frac{1}{2}\left(
[\omega(\al),[\omega(\beta),\omega(\gamma)]]-
[\omega([\al,\beta]),\omega(\gamma)]\right)\\
-\bigodot_{\al,\beta,\gamma}\frac{1}{2}\left(
[\omega(\al),\omega([\al,\beta])]+
\omega([[\al,\beta],\gamma])\right),
\end{multline*}
where the symbol $\bigodot$ denotes cyclic sum over the subscripts.
The first and fourth term vanish because of Jacobi's identity,
while the two middle terms cancel out.
\end{proof}

Let $s_j: U_j\to P$ be a local section of $P(M,G)$ defined over an
open set $U_j\subset M$. Then we have a trivializing isomorphism
$\psi_j:p^{-1}(U_j)\to U_j\times G$ such that
$s_j(x)=\psi^{-1}_j(x,e)$, where $e\in G$ is the identity. If $\gg$
denotes the Lie algebra of $G$, we also have an isomorphism
$T_{U_j}P/G\simeq TU_j\times\gg$. Given an $A$-connection $h$, the
corresponding connection 1-section $\omega$ trivializes over $U_j$
as
\[ \omega(\al)\simeq(\#\al,\omega_j(\al)),\]
for some $\omega_j\in\Gamma(A^*,\gg)$. We also have the following
alternative description of $\omega_j$: if $\al\in\Gamma(A)$, $x\in
U_j$, and $u=s_j(x)$, then
\[X_u=(s_j)_*\#\al_x-h(s_j(x),\al_x)\in T_u P\]
is a vertical vector since, by (CI), we have:
\[ p_*X_u=p_*\cdot(s_j)_*\#\al_x-p_* h(s_j(x),\al_x)=\#\al_x-\#\al_x=0.\]
Then $\omega_j(\al)_x$ is the unique element $B\in\gg$ such that
$X_u=\sigma(B)_u$ , which exists by (CII). The
$\set{\omega_j}$ will be called the
\emph{local connection 1-sections} of the $A$-connection.

If $s_k: U_k\to P$ is another local section with $U_j\cap
U_k\not=\emptyset$, we denote by $\psi_{jk}:U_j\cap U_k
\to G$ the corresponding transition function. The
following proposition gives the transformation rule for the
local connection 1-sections. The proof is similar to the proof for the
Poisson case (\cite{Fernandes:article:1}, Prop.~1.3.1) and so it
will be omitted.

\begin{proposition}
The local connection 1-sections $\set{\omega_j}$ transform by
\begin{equation}
\label{eq:transform:connection}
\omega_k=\text{Ad}(\psi_{jk}^{-1})\omega_j+\psi_{jk}^{-1}d_A\psi_{jk},
\quad \text{on } U_j\cap U_k.
\end{equation}
Conversely, given a family of $\gg$-valued 1-sections
$\set{\omega_j}$, each defined in $U_j$ and satisfying relations
(\ref{eq:transform:connection}), there is a unique $A$-connection
in $P(M,G)$ which gives rise to the $\set{\omega_j}$.
\end{proposition}

For the local description of the curvature we observe that the
2-section $\Omega$ is vertical: in fact, by
(\ref{eq:curvature:geometric}), for $\al,\beta\in\Gamma(A)$ we have
\begin{align*}
p_*\Omega(\al,\beta)&=
\frac{1}{2}\left([p_*\omega(\al),p_*\omega(\beta)]-p_*\omega([\al,\beta])\right)\\
&=\frac{1}{2}\left([\#\al,\#\beta]-\#[\al,\beta]\right)=0.
\end{align*}
Hence, over each trivializing neighborhood $U_j$, the curvature
2-section trivializes as
\[ \Omega(\al,\beta)\simeq(0,\Omega_j(\al,\beta)),\]
for some $\Omega_j\in\Gamma^2(A^*,\gg)$.

The \emph{local curvature 2-sections} $\set{\Omega_j}$ satisfy local
versions of the structure equation (\ref{eq:first:structure:equation})
and Bianchi identity (\ref{eq:Bianchi}). Again, the proof is similar
to the Poisson case (\cite{Fernandes:article:1}, Prop.~1.4.1) and will
be omitted.

\begin{proposition}
The local curvature 2-sections of a connection transform by
\begin{equation}
\label{eq:transform:curvature}
\Omega_k=\text{Ad}(\psi_{jk}^{-1})\Omega_j,\quad \text{on }
U_j\cap U_k.
\end{equation}
Moreover, they are related to the local 1-sections by the first
structure equation
\begin{equation}
\label{eq:first:structure:equation:1}
\Omega_j=d_A \omega_j+\frac{1}{2}[\omega_j,\omega_j].
\end{equation}
and they satisfy the Bianchi identity:
\begin{equation}
\label{eq:Bianchi:1}
d_A \Omega_j+[\omega_j,\Omega_j]=0.
\end{equation}
\end{proposition}

Note that since the curvature is a vertical 2-section, given
$\al,\beta\in\Gamma(A)$ we can identify
$\Omega(\al,\beta)$ with a $\gg$-valued map in $P$. Under this
identification relation (\ref{eq:curvature:geometric}) can be written
as:
\begin{equation}
[h(\al),h(\beta)]-h([\al,\beta])=-2\sigma(\Omega(\al,\beta)),
\end{equation}
and we have
\[ \Omega(\al,\beta)_{s_j(x)}=\Omega_j(\al,\beta)_x.\]
Later we shall use this identification without further notice.

\subsection{Parallelism and Holonomy}
If $\gamma:[0,1]\to M$ is a smooth curve lying on a leaf $L$ of the
Lie algebroid $A$, then $\gamma$ is also smooth as map
$\gamma:[0,1]\to L$. This follows from the existence of ``canonical
coordinates'' for $M$ given by the Local Splitting Theorem. Also, by the
same theorem, we can choose (not uniquely) a piecewise smooth
family $t\mapsto\al(t)\in A$ such that $\#\al(t)=\dot{\gamma}(t)$.
Recalling definition \ref{defn:A:path}, this means that any path
that lies on a leaf is the base path of some $A$-path. Clearly, if
$\#$ is not injective, different $A$-paths can have the same base
path.

Let $\al(t)$ be an $A$-path with base path $\gamma(t)$. For any
$u_0\in P$ with $p(u_0)=\gamma(0)$ one can show, using the
$G$-invariance of $h$, that there exists a unique horizontal lift
$\tilde{\gamma}:[0,1]\to P$, which satisfies the system
\begin{equation}
\label{eq:lift}
\left\{
\begin{array}{l}
\dot{\tilde{\gamma}}(t)=h(\tilde{\gamma}(t),\al(t)),\\
\\
\tilde{\gamma}(0)=u_0.
\end{array}
\right.
\end{equation}
Hence, we can define parallel displacement of the fibers along an
$A$-path $\al(t)$ in the usual way: if $u_0\in p^{-1}(\gamma(0))$
we define $\tau(u_0)=\tilde{\gamma}(1)$, where $\tilde{\gamma}(t)$
is the unique horizontal lift of $\al(t)$ starting at $u_0$. We
obtain a map $\tau:p^{-1}(\gamma(0))\to p^{-1}(\gamma(1))$, which
will be called \emph{parallel displacement} of the fibers along the
$A$-path $\al(t)$. It is clear, since horizontal curves are mapped
by $R_a$ to horizontal curves, that parallel displacement commutes
with the action of $G$:
\begin{equation}
\tau\circ R_a=R_a\circ\tau.
\end{equation}
Therefore, parallel displacement is an isomorphism between the
fibers.

If $x\in M$ belongs to the leaf $L$, let $\Omega(L,x)$ be the loop
space of $L$ at $x$. An $A$-path $\al(t)$ for which the base path is a
loop $\gamma\in\Omega(L,x)$ will be called an \emph{$A$-loop} in
$L$ based at $x$. Parallel displacement along such an $A$-loop
$\al(t)$ gives a an isomorphism of the fiber $p^{-1}(x)$ into
itself. The set of all such isomorphisms forms the holonomy group
of the $A$-connection, with reference point $x$, and is denoted
$\Phi(x)$. Similarly, one has the restricted holonomy group, with
reference point $x$, denoted $\Phi^0(x)$, which is defined by using
$A$-loops in $L$ with base paths homotopic to the constant path.

For any $u\in p^{-1}(x)$ we can also define the holonomy groups
$\Phi(u)$ and $\Phi^0(u)$. Just as in the covariant case, $\Phi(u)$
is the subgroup of $G$ consisting of those elements $a\in G$ such
that $u$ and $ua$ can be joined by an horizontal curve. We have
that $\Phi(u)$ is a Lie subgroup of $G$, with connected component
of the identity $\Phi^0(u)$, and we have isomorphisms
$\Phi(u)\simeq \Phi(x)$ and $\Phi(u)^0\simeq\Phi(x)^0$.

If $x,y\in M$ belong to the same leaf then the holonomy groups
$\Phi(x)$ and $\Phi(y)$ are isomorphic. This is because if $u,v\in
P$ are points such that, for some $a\in G$, there exists an
horizontal curve connecting $ua$ and $v$, then
$\Phi(v)=Ad(a^{-1})\Phi(u)$, so $\Phi(u)$ and $\Phi(v)$ are
conjugate in $G$. However, if $x,y\in M$
belong to different leaves the holonomy groups $\Phi(x)$ and $\Phi(y)$
will be, in general, non-isomorphic.

The holonomy groups can be given an infinitesimal description as in
the Ambrose-Singer Holonomy Theorem. For that, suppose that
$\gamma\in A_{x}$ satisfies $\#\gamma=0$. If $u\in p^{-1}(x)$ we
set:
\[ \Lambda(\gamma)_{u}\equiv \omega_j(\gamma)_x,\]
where $j$ is such that $s_j(x)=u$. It follows from the
transformation rule (\ref{eq:transform:connection}) 
that the previous formula gives a well defined map $\Lambda: \Omega(L, x)
\to \mathfrak{g}$ (the dashed arrow in the diagram above). Also, denote
by $P(u_0)$ the set of points $u\in P$ that can be joined to $u_0$
by an horizontal curve. We have:

\begin{theorem}[Holonomy Theorem]
\label{thm:holonomy}
Given any $A$-connection in the principal bundle
$P(M,G)$ and $u_0\in P$, the Lie algebra of the holonomy group
$\Phi(u_0)\subset G$ is the ideal of $\gg$ spanned by all elements 
$\Omega(\al,\beta)_{u}$ and $\Lambda(\gamma)_{u}$, 
where $u\in P(u_0)$ and $\al,\beta,\gamma\in A_{p(u)}$, with $\#\gamma=0$.
\end{theorem}

The proof is also analogous to the Poisson case
(\cite{Fernandes:article:1}, Thm.~1.5.2) and so it will be omitted.

Note the presence of the extra term $\Lambda$ in the Holonomy
Theorem. This means that a Lie algebroid connection can be flat and
still have non-discrete holonomy, a phenomenon that is not present
in the covariant case or whenever $\#$ is injective.

\subsection{Relationship to Ordinary Connections}
Consider the tangent Lie algebroid $A=TM$ and a
$TM$-co\-nnec\-tion in $P(M,G)$, i.e., a
covariant connection. Its horizontal lift $\bar{h}:p^*TM\to TP$
is completely determined by the horizontal distribution $\Ho$.
For a Lie algebroid $A\to M$, the formula $h(u,\al)\equiv \bar{h}(u,\#\al)$
defines an $A$-connection in $P(M,G)$ which is said to be \emph{induced
by the covariant connection}. Note that in this case the lift $h$
satisfies:
\begin{equation}
\label{eq:F:connection}
\#\al=0\Longrightarrow h(u,\al)=0, \qquad (u,\al)\in p^*A.
\end{equation}
This construction shows that there are always $A$-connections on any
principal bundle $P(M,G)$ over a Lie algebroid $A\to M$.

Let $\bar{\omega}$ be the connection 1-form and let
$\bar{\Omega}$ be the curvature 2-form of the covariant connection $\bar{h}$.
Then it is clear from the definitions given above that the
connection 1-section $\omega$ and the curvature
2-section $\Omega$ of the induced $A$-connection $h$ are given by:
\begin{equation}
\omega=\#^*\bar{\omega},\qquad \Omega=\#^*\bar{\Omega}.
\end{equation}
Also, given trivialization isomorphisms $\set{\psi_j}$, inducing local sections
$\set{s_j}$, we see that the associated local connection 1-sections
and curvature 2-sections are related by:
\begin{equation}
\omega_j=\#^*\bar{\omega}_j,\qquad \Omega_j=\#^*\bar{\Omega}_j.
\end{equation}

However, in general, a connection will not satisfy property
(\ref{eq:F:connection}) and we set:

\begin{definition}
An $A$-connection on a principal bundle $P(M,G)$
is called an \textsc{$\F$-connection} if its horizontal lift satisfies
condition (\ref{eq:F:connection})
\end{definition}

Let us fix one such $\F$-connection on $P(M,G)$. Then, on
the pull-back bundle $p_L:i^*P\to L$, we have a naturally induced
$TL$-connection, i.e., a covariant connection.
On the total space 
\[ i^*P= \set{(x,u)\in
L\times P: i(x)=p(u)}\] 
we define the horizontal lift
$\bar{h}_L:p_L^*TL\to T(i^*P)$ by setting
\begin{equation}
\label{eq:S:connection}
\bar{h}_L((x,u),v)=(v,h(u,\al)), \qquad (x,u)\in i^*P,
v\in T_x L,
\end{equation}
where we choose any $\al\in A_x$ such that $\#\beta=v$, and
we are identifying $T(i^*P)=\set{(v,w)\in TS\times TP:v=p_*w}$. Note
that if $\#\beta'=\#\beta=v$ we get the same result in
(\ref{eq:S:connection}) since $h$ is an $\F$-connection.

\begin{proposition}
\label{prop:F:connections}
Let $\Omega$ and $\omega$ be the connection and curvature sections
for an $\F$-connection in $P(M,G)$. For a leaf $i:L\hookrightarrow M$
denote by $\bar{\omega}^L$ and $\bar{\Omega}^L$ the connection
1-form and the curvature 2-form for the induced connection on
$i^*P(M,G)$. Then $\omega$ and $\Omega$ are
$i$-related to $\#^*\bar{\omega}^L$ and $\#^*\bar{\Omega}^L$:
\begin{equation}
\label{eq:relationship:connection}
i_*\#^*\bar{\omega}^L=\omega,\qquad i_*\#^*\bar{\Omega}^L=\Omega.
\end{equation}
\end{proposition}

Therefore, an $\F$-connection can be thought of as a \emph{family}
of ordinary connections over the leaves of $M$. The
connection 1-section $\omega$ and the curvature
2-sections $\Omega$ are obtained by gluing together the
connection 1-sections $\#^*\bar{\omega}^L$ and the curvature 2-sections
$\#^*\bar{\Omega}^L$ of the connections on the leaves of $M$.

For an $\F$-connection, horizontal lifts of $A$-paths
$\al(t)$ depend only on the base path $\gamma(t)$. Therefore, one
has a well determined notion of horizontal lift of a curve lying on a leaf.
It follows that for these connections, parallel
displacement can also be defined by first reducing to the pull-back
bundle over a leaf and then parallel displace the
fibers. Hence, the holonomy groups $\Phi(x)$ and $\Phi^0(x)$
coincide with the usual holonomy groups of the pull-back connection
on the leaf $L$ through $x$.

\subsection{Connections on a Vector Bundle}

If $G$ acts on the left on a manifold $F$ we denote by
$p_E:E(M,F,G,P)\to M$ the fiber bundle associated with $P(M,G)$
with standard fiber $F$.

Given an $A$-connection in $P(M,G)$ with horizontal
lift $h:p^*A\to TP$, we define the induced horizontal lift
$h_E:p_E^*A\to TE$ as follows: given $w\in E$ choose $(u,\xi)\in
P\times F$ which is mapped to $w$, and set
\begin{equation}
\label{eq:connection:principal:fiber}
h_E(w,\al)\equiv\xi_*h(u,\al),
\end{equation}
where we are identifying $\xi$ with the map $P\to E$ which sends an
element $u\in P$ to the equivalence class $[u,\xi]\in E$. One can
check easily that this definition does not depend on the choice of
$(u,\xi)$, so we obtain a well defined bundle map $h_E:p_E^*A\to TE$
which makes the following diagram commute:
\begin{equation}
\label{diagram:connection:fiber:space}
\xymatrix{
p_E^*A\ar[r]^{h_E}\ar[d]_{\widehat{p}_E}& TE \ar[d]^{p_{E*}} \\
A\ar[r]_{\#} & TM
}
\end{equation}

As before, we can define horizontal and vertical vectors in $TE$,
horizontal lifts to $E$ of curves lying on leaves of the orbit 
foliation, and parallel displacement of fibers of $E$. We shall call a
cross section $\sigma$ of $E$ over an open set $U\subset M$
\emph{parallel} if $\sigma_*(v)$ is horizontal for all tangent
vectors $v\in T_U M$.

Now assume that $G$ acts linearly on a vector space $V$. On the
associated vector bundle $E(M,V,G,P)$ we obtain an horizontal lift
$h_E:p_E^*A\to TE$ which has the distinguish property of being
\emph{linear}:
\begin{equation}
\label{eq:connection:linear}
h_E(w,a_1\al_1+a_2\al_2)=a_1h_E(w,\al_1)+a_2h_E(w,\al_2).
\end{equation}
Conversely, given a bundle map $h_E:p_E^*A\to TE$ satisfying
(\ref{eq:connection:linear}) and making the diagram
(\ref{diagram:connection:fiber:space}) commute, it is associated with
some $A$-connection on the principal bundle $P(M,G)$ through relation
(\ref{eq:connection:principal:fiber}).

For an $A$-connection on a vector bundle $E$ we have the notion of
\emph{$A$-derivative} of sections of $E$ along $A$-paths, analogous
to the notion of covariant derivative of sections for covariant connections.
Given a section $\phi$ of $E$ defined along an $A$-path $\al(t)$, its
\emph{$A$-derivative} $\nabla_{\al}\phi$ is the section defined by
\begin{equation}
\label{eq:contra:derivative}
\nabla_{\al}\phi(t)\equiv\lim_{h\to 0}\frac{1}{h}
\left[\tau^{t+h}_t(\phi(\gamma(t+h)))-\phi(\gamma(t))\right],
\end{equation}
where $\gamma(t)$ denotes the base path of $\al(t)$ and
$\tau^{t+h}_t:p_E^{-1}(\gamma(t+h))\to p_E^{-1}(\gamma(t))$
denotes parallel transport of the fibers from $\gamma(t+h)$ to
$\gamma(t)$ along the $A$-path. The proof of the
following proposition is elementary.

\begin{proposition}
\label{prop:derivative:properties}
Let $\phi$ and $\psi$ be sections of $E$ and $f$ a function on $M$
defined along $\gamma$. Then
\begin{enumerate}
\item[(i)] $\nabla_{\al}(\phi+\psi)=\nabla_{\al}\phi+\nabla_{\al}\psi$;
\item[(ii)] $\nabla_{\al}(f\phi)=(f\circ\gamma)\nabla_{\al}\phi+
\dot{\gamma}(f)(\phi\circ\gamma)$;
\end{enumerate}
\end{proposition}

Now let $\al\in A_x$ and let $\phi$ be a cross section of
$E$ defined in a neighborhood of $x$. The $A$-derivative
$\nabla_\al\phi$ of $\phi$ in the direction of $\al$ is defined as
follows: choose an $A$-path $\al(t)$, with base path $\gamma(t)$, defined for
$t\in(-\eps,\eps)$ and such that $\gamma(0)=x$, $\al(0)=\al$.
Then we set:
\begin{equation}
\label{eq:contra:derivative:0}
\nabla_\al\phi\equiv\nabla_{\al(t)}\phi(0).
\end{equation}
It is easy to see that $\nabla_\al\phi$ is independent of the choice of
$A$-path. Clearly, a cross section $\phi$ of $E$ defined on an open
set $U\subset M$ is flat iff $\nabla_\al\phi=0$ for all $\al\in A_x$, $x\in
M$.

Finally,  given $\al\in\Gamma(A)$ and $\phi$ a section of
$E$,  we can also define the $A$-derivative $\nabla_\al\phi$ to be the
section of $E$ given by:
\begin{equation}
\label{eq:contra:derivative:00}
\nabla_\al\phi(x)=\nabla_{\al_x}\phi.
\end{equation}
Moreover, we have the following properties of the $A$-derivative:

\begin{proposition}
\label{prop:derivative:properties:2}
The $A$-derivative $\nabla$ is a map $\Gamma(A)\times\Gamma(E)\to \Gamma(E)$ that
satisfies
\begin{enumerate}
\item[(i)] $\nabla_{\al+\beta}\phi=\nabla_\al\phi+D_\beta\phi$;
\item[(ii)] $\nabla_\al(\phi+\psi)=\nabla_\al\phi+\nabla_\al\psi$;
\item[(iii)] $\nabla_{f\al}=f\nabla_\al\phi$;
\item[(iv)] $\nabla_\al(f\phi)=f\nabla_\al\phi+\#\al(f)\phi$;
\end{enumerate}
for all $\al,\beta\in\Gamma(A)$, $\phi,\psi\in\Gamma(E)$, and $f\in
C^{\infty}(M)$.
\end{proposition}

It is also true that the $A$-derivative uniquely determines the
connection: Given a map $\nabla:\Gamma(A)\times\Gamma(E)\to\Gamma(E)$
satisfying properties (i) to (iv) of 
Proposition \ref{prop:derivative:properties:2},
there exists a unique $A$-connection on the associated principal
bundle $P(M,G)$ for which the induced $A$-derivative on $E$ is $\nabla$.

For an $A$-connection in a vector bundle $E$ we define the
\emph{curvature section} $R$ to be the section of
$\wedge^2 A^*\otimes\text{End}(E)$ given by
\begin{equation}
\label{curvature}
R(\al,\beta)\gamma\equiv s_j(x)\left[\Omega_j(\al,\beta)\cdot
s_j^{-1}(x)(\gamma)\right].
\end{equation}
where $x\in U_j$, $\al,\beta\in A_x$ and $\gamma\in E_x$ (here we view $u\in
P(M,G)$ as an isomorphism $u: V\to E_{p(u)}$).
Note that if $x\in U_j\cap U_k$ and
$s_k(x)=\psi_{jk}(x)s_k(x)$ we obtain the same values in formula
(\ref{curvature}), so this really defines
a global section on $M$. This section can be easily
expressed in terms of $A$-derivatives as
\begin{equation}
\label{eq:curvature:derivative}
R(\al,\beta)\gamma=
\nabla_\al \nabla_\beta \gamma-\nabla_\beta \nabla_\al \gamma
-\nabla_{[\al,\beta]}\gamma.
\end{equation}
Moreover, Bianchi's identity (\ref{eq:Bianchi:1}) in this notation
reads
\begin{equation}
\label{eq:Bianchi:1'}
\bigodot_{\al_1,\al_2,\al_3} \nabla_{\al_1}(R(\al_2,\al_3))-
\bigodot_{\al_1,\al_2,\al_3} R([\al_1,\al_2],\al_3)=0.
\end{equation}

If the $A$-connection $h$ is induced by a
covariant connection $\bar{h}$, the $A$-derivative $\nabla$ and the
covariant derivative $\bar{\nabla}$ are related by
\begin{equation}
\nabla_\al =\bar{\nabla}_{\#\al}.
\end{equation}
On the other hand, $\F$-connections can be characterized by the
condition:
\begin{equation}
\#\al=0\ \Longrightarrow\ \nabla_\al=0, \qquad \forall \al\in \Gamma(A).
\end{equation}
Moreover, by Proposition \ref{prop:F:connections}, for an
$\F$-connection, on each leaf $i:L\hookrightarrow M$
there is a covariant connection on the pullback bundle $i^*P$,
inducing a covariant derivative $\nabla^L$ on $i^*E$, with the
following property: if $\psi$ is any cross section of $E$, then
\begin{equation}
i^*\nabla_\al \psi =\nabla^L_{\#i^*\al}i^*\psi,
\end{equation}
where $i^*\psi$ denotes the section of the pullback bundle $i^*E$
induced by $\psi$.

\begin{demo}{Remark}
A flat $A$-connection on a vector bundle $E$
is sometimes called a \emph{Lie algebroid representation} of $A$ or an
$A$-module (see e.g.~\cite{Evens:article:1,Huebschmann:article:2}).
In fact, if we set $\al\cdot
s\equiv\nabla_\al s$ we get a bilinear product
$\Gamma(A)\times\Gamma(E)\to\Gamma(E)$ and the axioms for
$\nabla$ to be a flat $A$-connection are translated into
\begin{align}
(f\al)\cdot s&=f(\al\cdot s),\\
\al\cdot(fs)&=f(\al\cdot s)+\#\al(f)s,\\
[\al,\beta]\cdot s&=\al\cdot(\beta\cdot s)-\beta\cdot(\al\cdot s).
\end{align}
One reason for this terminology is that in the case of a Lie algebra
$\gg$, viewed as a Lie algebroid over a one point space, $E$ is just a
vector space and these are the requirements for $E$ to be a $\gg$-module.
\end{demo}

\subsection{Linear $A$-Connections}
\label{section:linear:connections}

A \emph{linear $A$-connection} is a $A$-connection on the
frame bundle $P=GL(A)$, so $G=GL(r)$ where $r=\rank A$.
If $u=(\al_1,\dots,\al_r)\in GL(A)$ is a frame, we can view $u$ as a
linear isomorphism $u:\Rr^r \to A_{p(u)}$ by setting
\[ u(v_1,\dots,v_r)=\sum_{s}v_s\al_s, \qquad (v_1,\dots,v_r)\in\Rr^r.\]
We define the \emph{canonical 1-sections} $\theta_j\in \Gamma(A^*)$ on
an open set $U_j$, with trivializing isomorphism
$\psi_j:p^{-1}(U_j)\to U_j\times G$, and associated section
$s_j(x)=\psi_j^{-1}(x,e)$, to be the $\Rr^r$-valued 1-sections defined
by
\begin{equation}
\label{eq:canonical:fields}
\theta_j(\al)_x=s_j(x)^{-1}(\al), \qquad x\in U_j.
\end{equation}
Given an $A$-connection these allow us to define the \emph{torsion 2-sections}
$\Theta_j\in\Gamma(\wedge^2 A^*)$ to be the $\Rr^r$-valued 2-sections
given by
\begin{equation}
\label{eq:torsion:fields}
\Theta_j(\al,\beta)=d_A\theta_j(\al,\beta)+
\omega_j(\al)\cdot\theta(\beta)-\omega_j(\beta)\cdot\theta_j(\al).
\end{equation}

We state the main properties of the torsion:
\begin{proposition}
The canonical 1-sections and the torsion 2-sections of a linear
$A$-connection are related by
\begin{eqnarray}
\label{eq:transform:canonical}
  \theta_k=\psi_{jk}^{-1}\cdot\theta_j,\\
\label{eq:transform:torsion}
  \Theta_k=\psi_{jk}^{-1}\cdot\Theta_j.
\end{eqnarray}
Moreover, they satisfy the Bianchi identity
\begin{equation}
\label{eq:Bianchi:2}
d_A \Theta_j(\al,\beta,\gamma)=
\bigodot_{\al,\beta,\gamma}d_A\omega_j(\al,\beta)\cdot\theta_j(\gamma)-
\bigodot_{\al,\beta,\gamma}\omega_j(\al)\cdot d_A\theta_j(\beta,\gamma).
\end{equation}
\end{proposition}

The vector bundle $A$ is associated with the principal bundle
$GL(A)$ of frames of $A$. Therefore, as it was explained in the
previous paragraph, any linear
$A$-connection determines an $A$-derivative operator
$\nabla:\Gamma(A)\times\Gamma(A)\to\Gamma(A)$ such that:
\begin{eqnarray}
\label{eq:contra:derivative:1}
\nabla_{f_1\al_1+f_2\al_2}=f_1 \nabla_{\al_1}+f_1 \nabla_{\al_1}, \qquad
\text{for all } f_i\in C^\infty(M),\ \al_i\in\Gamma(A),\\
\label{eq:contra:derivative:2}
\nabla_\al (f\beta)=f \nabla_\al \beta+\#\al(f)\beta,\qquad \text{for all }
f\in C^\infty(M),\
\al,\beta\in\Gamma(A).
\end{eqnarray}

One can also consider other associated vector bundles to $GL(A)$
which lead, just us in the covariant case, to $A$-derivatives of
any r-sections over $A$. For example, if $X$ is a section of $A^*$,
then $\nabla_\al X$, the \emph{$A$-derivative of $X$} along
$\al\in\Gamma(A)$, is completely characterized by the relation
\begin{equation}
\seq{\nabla_\al X,\beta}=\#\al(\seq{X,\beta})-\seq{X,\nabla_\al \beta},
\end{equation}
which must hold for every section $\beta\in\Gamma(A)$.

For a linear $A$-connection we define the \emph{torsion section} $T$
to be the section of $A\otimes A^*\otimes A^*$ given by
\begin{equation}
\label{torsion}
T(\al,\beta)\equiv s_j(x)(\Theta_j(\al,\beta)),
\end{equation}
where $x\in U_j$, $\al,\beta,\gamma\in A_x$. Note that if
$x\in U_j\cap U_k$ and
$s_k(x)=\psi_{jk}(x)s_k(x)$ we obtain the same values in formula
(\ref{torsion}), so this really defines
a global section of $M$. In terms of $A$-derivatives, the torsion is
given by
\begin{equation}
\label{eq:torsion:derivative}
T(\al,\beta)=\nabla_\al \beta-\nabla_\beta \al - [\al,\beta].
\end{equation}
Moreover, the Bianchi identities
(\ref{eq:Bianchi:1'}) and (\ref{eq:Bianchi:2}) can also be expressed as
\begin{eqnarray}
\bigodot_{\al,\beta,\gamma}\left(\nabla_\al R(\beta,\gamma)+
R(T(\al,\beta),\gamma)\right)=0,\\
\bigodot_{\al,\beta,\gamma}\left(R(\al,\beta)\gamma-
T(T(\al,\beta),\gamma)-\nabla_\al T(\beta,\gamma)\right)=0.
\end{eqnarray}

If it happens that the $A$-connection is related to some
covariant connection by:
\[ \#\nabla_\al\beta =\bar{\nabla}_{\#\al}\#\beta,\]
the torsion and curvature sections of the $A$-connections
are transformed by the musical
homomorphism to the usual torsion and tensor fields of $\nabla$:
\begin{align*}
\bar{T}(\#\al,\#\beta)&=\#T(\al,\beta),\\
\bar{R}(\#\al,\#\beta)\#\gamma&=\#R(\al,\beta)\gamma.
\end{align*}

Local coordinate expressions for linear $A$-connections can be
obtained in a way similar to the covariant case. Let $(x^1,\dots,x^m)$
be local coordinates and $(\al^1,\dots,\al^r)$ local sections for $A$
in a trivializing neighborhood $U\subset M$. Then we can define
\emph{Christoffel symbols} $\Gamma^{st}_u$ by
\begin{equation}
\nabla_{\al^s}\al^t =\Gamma^{st}_u \al^u,
\end{equation}
where we are using the repeated index sum convention.  It is easy to see
that under a change of coordinates and a change of basis of the form
\[ y^j=y^j(x^1,\dots,x^m),\qquad
\tilde{\al}^{s'}=a^{s'}_s(x^1,\dots,x^m)\al^s,\]
these symbols transform according to
\begin{equation}
\label{eq:symbols:transform}
{\tilde{\Gamma}}^{s' t'}_{u'}=
a^{s'}_s a^{t'}_t \tilde{a}^{u}_{u'}\Gamma^{st}_u
+
a^{s'}_s b^{s i}\frac{\partial a^{t'}_t}{\partial x^i}\tilde{a}^{t}_{u'}
\end{equation}
where $b^{s i}$ are the structure functions of $\#$ in the original
coordinates, and $(\tilde{a}^{s}_{s'})$ denotes the inverse of
$(a^{s'}_s)$.  Conversely, given a family of symbols that transform
according to this rule under a change of coordinates and basis of
sections, we obtain a well defined $A$-derivative/connection on $A$.

Let us call an element
\[ K\in \Gamma(\underbrace{A\otimes\dots\otimes A}_{k\text{ times}}
\otimes\underbrace{A^*\otimes\dots\otimes A^*}_{l\text{ times}})\]
a section of type $(k,l)$, or simply a $(k,l)$-section. Using these
symbols, it is easy to get the local coordinates expressions for the
$A$-derivatives of a such a $(k,l)$-section: Let $\al=\sum_s c_s
\al^s$ be a section of $A$ and write $K$ as
\[ K=K^{t_1\dots t_l}_{s_1\dots s_k}
\al^{s_1}\otimes\cdots\otimes\al^{s_k}
\otimes
\check{\al}_{t_1}\otimes\cdots\otimes\check{\al}_{t_l}\]
where $\set{\check{\al}_{1},\dots,\check{\al}_{r}}$ is the basis of
local sections of $A^*$ dual to $\set{\al^1,\dots\al^r}$. Then the
$A$-derivative of $K$ along $\al$ is the $(k,l)$-section
\begin{multline}
\label{eq:derivative}
(\nabla_\al K)^{s_1\dots s_l}_{t_1\dots t_k}=
b^{sj}\al_s\frac{\partial K^{s_1\dots s_l}_{t_1\dots t_k}}{\partial x^j}
-\sum_{a=1}^l\left(\Gamma^{s s_a}_u\al_s K^{s_1\dots u\dots
s_l}_{t_1\dots t_s}\right)\\ +\sum_{b=1}^k\left(\Gamma^{s
u}_{t_b}\al_s K^{s_1\dots s_l}_{t_1\dots u\dots t_s}\right).
\end{multline}

Given a section $K$ of type $(k,l)$ we shall write, as in the
covariant case, $\nabla K$ for the unique section of type $(k,l+1)$
that satisfies
\begin{equation}
(\nabla K)^{s_1\dots s_l, s}_{t_1\dots t_k}=
(\nabla_{\al^s} K)^{s_1\dots s_l}_{t_1\dots t_k}.
\end{equation}
A tensor field $K$ on $M$ is \emph{parallel} iff $\nabla K=0$.

From formulas (\ref{eq:torsion:derivative}) and
(\ref{eq:curvature:derivative}), we obtain immediately the following
local coordinates expressions for the torsion and curvature
in terms of Christoffel symbols and structure functions:
\begin{align}
\label{eq:torsion:coordinates}
T^{st}_u&=\Gamma^{st}_u-\Gamma^{st}_u-c^{st}_u,\\
\label{eq:curvature:coordinates}
R^{stu}_v&=\Gamma^{sa}_v\Gamma^{tu}_a-\Gamma^{ta}_v\Gamma^{su}_a
+b^{si}\frac{\partial\Gamma^{tu}_v}{\partial
x^i}-b^{ti}\frac{\partial\Gamma^{su}_v}{\partial x^i}
-c^{st}_a\Gamma^{au}_v.
\end{align}

\subsection{Connections Compatible with the Lie Algebroid Structure}
\label{sec:compatible}
In Poisson geometry, linear connections for which the Poisson tensor
is parallel play an important role. We recall that for a Poisson manifold
$(M,\Pi)$ a contravariant connection is just a linear connection
$\nabla$ on $T^*M$. It induces a $T^*M$-connection
$\check{\nabla}$ on $TM$ in the usual way:
\[ \langle \check{\nabla}_\al X,\beta\rangle=
\#\al(\seq{X,\beta})-\seq{X,\nabla_\al \beta},\qquad
\forall X\in \X^1(M),\al,\beta\in\Omega^1(M).\]
The connection $\nabla$ is a Poisson connection, in the sense that
$\nabla\Pi=0$, iff
\[ \check{\nabla}\#=\#\nabla.\]

For a general Lie algebroid we do not have a relationship between
$A$ and $TM$ analogue to the duality between $T^*M$ and $TM$ that
occurs in the Poisson case. The notion of a Poisson connection is
replaced by the following:

\begin{definition}
A linear connection on a Lie algebroid $A$, with associated
$A$-derivative $\nabla$, is said to be
\textsc{compatible with the Lie algebroid
structure} of $A$ if there exists an $A$-connection in $TM$, with
associated $A$-derivative $\check{\nabla}$, such that
\[ \check{\nabla}\#=\#\nabla.\]
\end{definition}

In \cite{Fernandes:article:1}, Prop.~2.5.1, a simple argument shows
that Poisson manifolds always admit Poisson connections.
This argument extends to Lie algebroids.

\begin{proposition}
\label{prop:compatible:connection}
Every Lie algebroid admits compatible linear connections.
\end{proposition}

\begin{proof}
Let $U_a$ be a domain of a chart $(x^1,\dots,x^m)$ where there
exists a basis of trivializing sections $\set{\al^1,\dots,\al^r}$.
On $U_a$, we define a linear $A$-connection by
\[ \nabla^{(a)}_{\al^s}\al^t=c^{st}_u\al^u,\]
and an $A$-connection on $TU_a$ by
\[ {\check{\nabla}}^{(a)}_{\al^s}\frac{\partial}{\partial x^i}
=\frac{\partial b^{sj}}{\partial x^i}\frac{\partial }{\partial x^j},\]
where $b^{sj}$ and $c^{st}_u$ denote, as usual, the structure
functions for this choice of coordinates and basis. A straight
forward computation shows that the relation
$\#[\al^s,\al^t]=[\#\al^s,\#\al^t]$ implies that
\[ {\check{\nabla}}^{(a)}\#=\#\nabla^{(a)},\]
so $\nabla^{(a)}$ is a linear connection in $U_a$ compatible with the Lie
algebroid structure.

If we take an open cover of $M$ by such chart domains and if
$\sum_a\phi^{(a)}=1$ is a partition of unity subordinated to this
cover, then $\nabla\equiv\sum_a \phi^{(a)} \nabla^{(a)}$ and
$\check{\nabla}\equiv\sum_a \phi^{(a)} {\check{\nabla}}^{(a)}$ define
$A$-connections that satisfy $\check{\nabla}\#=\#\nabla$,
i.e., $\nabla$ is a connection in $M$ compatible with the Lie
algebroid structure.
\end{proof}

An alternative approach to compatible connections is the following. Let
us write $E=A\ominus T^*M$ for the vector bundle over $M$ which is
the ``formal difference'' of the vector bundles $A$ and $T^*M$. If
we are given a linear $A$-connection compatible with the Lie
algebroid structure, we obtain an $A$-connection on $E$ by setting:
\begin{equation}
\tilde{\nabla}_\gamma (\al+\omega)\equiv \nabla_\gamma\al+
\check{\nabla}_\gamma\omega,
\qquad \gamma,\al\in \Gamma(A),\omega\in\Omega^1(M).
\end{equation}
Also on $E$ we have a \emph{canonical skew-symmetric bilinear form}
$(\cdot,\cdot)$ given by the formula:
\begin{equation}
\label{eq:bilinear:form}
(\al+\omega,\beta+\eta)\equiv \eta(\#\al)-\omega(\#\beta),
\end{equation}
and which, in general, will be degenerate.
\begin{proposition}
\label{prop:compatible}
$\nabla$ is compatible with the Lie algebroid structure iff
the skew-symmetric form $(\cdot,\cdot)$ is parallel with respect to
$\tilde{\nabla}$, i.e.,
\[\#\gamma((\zeta,\xi))=(\tilde{\nabla}_\gamma\zeta,\xi)+
(\zeta,\tilde{\nabla}_\gamma\xi),
\qquad \gamma\in\Gamma(A),\ \zeta,\xi\in\Gamma(E).\]
\end{proposition}

\begin{proof} A straightforward computation.
\end{proof}

\section{Holonomy}

For a regular foliation the topological behaviour close to a given
leaf is controlled by the holonomy of the leaf. For singular
foliations the situation is more complex (see e.g.~\cite{Dazord:article:1},
where holonomy is defined for transversally stable leaves). 

In this section, we will show that for any Lie algebroid it is
possible to introduce a notion of holonomy. This holonomy can be
defined as a map between the transversal algebroid germs that describe
the transversal geometry of the Lie algebroid (cf.~Section
\ref{sec:trasnverse:structure}). In this theory of holonomy
$A$-connections play a crucial role.

Later in the section we consider linear holonomy which will take us to
the concept of a basic connection. Basic connections will be used in
the next section to define secondary characteristic classes for the
orbit foliation of a Lie algebroid.

\subsection{Holonomy of a Leaf} 
Throughout this discussion we will consider a fixed leaf
$i:L\hookrightarrow M$ of the Lie algebroid $\pi:A\to M$.
We denote by $\nu(L)=T_LM/TL$ the normal
bundle to $L$ and by $p:\nu(L)\to L$ the natural projection.
By the Tubular Neighborhood Theorem, there exists a smooth immersion
$\tilde{i}:\nu(L)\to M$ satisfying the following properties:
\begin{enumerate}
\item[(i)] $\tilde{i}|_Z=i$, where we identify the zero section $Z$ of
$\nu(L)$ with $L$;
\item[(ii)] $\tilde{i}$ maps the fibers of $\nu(L)$ transversally to the
foliation of $M$;
\end{enumerate}

Assume that we have fixed such an immersion, and let $x\in L$. Each
fiber $F_x=p^{-1}(x)$ is a submanifold of $M$ transverse to the
foliation so we have (see Section \ref{sec:trasnverse:structure})
the transverse Lie algebroid structure $A_{F_x}\to F_x$. Because
$F_x$ is a linear space we can choose a trivialization and identify
the fibers $(A_{F_x})_u$ for different $u\in p^{-1}(x)$. Finally,
we choose a complementary vector subbundle $E\subset A$ to
$A_{F_x}$:
\begin{equation}
\label{eq:decompose}
A_u=E_u\oplus (A_{F_x})_u.
\end{equation}
Note that, by construction, the anchor $\#:A\to TM$ maps
$A_{F_x}$ onto $TF_x$, its restriction to $E$ is injective, and vectors
in $\#E$ are tangent to the orbit foliation.

Let $\al\in A_x$. We decompose $\al$ according to
(\ref{eq:decompose}):
\[\al=\al^\parallel+\al^{\perp},\text{where } \al^\parallel\in E_x,
\quad \al^{\perp}\in (A_{F_x})_x.\]
For each $u\in F_x=p^{-1}(x)$, we denote by
$\tilde{\al}^\parallel_u\in E_u$ the unique element such that $d_u
p\cdot \#\tilde{\al}^\parallel_u=\#\al^\parallel$, and by
$\tilde{\al}^\perp_u\in (A_{F_x})_u$ the element corresponding to
$\al^{\perp}$ under the identification $(A_{F_x})_u\simeq
(A_{F_x})_x$. We also set $\tilde{\al}\equiv
\tilde{\al}^\parallel+\tilde{\al}^\perp$.

Given $\al\in A_x$, $x\in L$, and $u\in F_x$, define the
\emph{horizontal lift} to $\nu(L)$ by
\[ h(u,\al)=\#\tilde{\al}_u\in T_u\nu(L).\]
By construction, we have the defining property of an $A$-connection:
\[ p_*h(u,\al)=\#\al,\quad u\in p^{-1}(x).\]
Note that $h$ depends on several choices made: tubular
neighborhood, trivialization of $A_{F_x}$ and complementary vector
bundle $E$.

Let $\al(t)$, $t\in [0,1]$, be an $A$-path with base path $\gamma(t)$
lying in the leaf $L$. If $u\in F_{\gamma(0)}=\nu(L)|_{\gamma(0)}$ is
a point in the fiber over $\gamma(0)$, there exists an $\eps>0$ and a
horizontal curve $\tilde{\gamma}(t)$ in $\nu(L)$, defined for
$t\in[0,\eps)$, which satisfies:
\[
\left\{
\begin{array}{l}
\frac{d}{dt}\tilde{\gamma}(t)=h(\tilde{\gamma}(t),\al(t)), \qquad
t\in[0,\eps),\\ \\
\tilde{\gamma}(0)=u.
\end{array}
\right.\]
If we take $u=0$ the lift $\tilde{\gamma}(t)$ coincides with
$\gamma(t)$, and so is defined for $t\in[0,1]$. It follows that we
can choose a neighborhood $U_\gamma$ of $0$ in
$F_{\gamma(0)}=\nu(L)|_{\gamma(0)}$, such that for each $u\in
U_\gamma$ the lift $\tilde{\gamma}(t)$ with initial point $u$ is
defined for all $t\in[0,1]$. Moreover, by passing from initial to
end point, this lift gives a diffeomorphism $H_L(\al)_0$ of
$U_\gamma$ onto a neighborhood $V_\gamma$ of $0$ in
$F_{\gamma(1)}=\nu(L)|_{\gamma(1)}$, with the property that $0$ is
mapped to $0$.

\begin{proposition}
\label{prop:holonomy}
Let $\al$ be an $A$-path in a leaf $L\subset M$.  The isomorphism
$H_L(\al)_0$ is covered by (a germ of) a Lie algebroid isomorphism
$H_L(\al)$ from $A_{F_{\gamma(0)}}$ to $A_{F_{\gamma(1)}}$.  If $\al'$
is another $A$-path in $L$ such that $\gamma(1)=\gamma'(0)$ we have
\begin{equation}
\label{eq:holonomy:homomorphism}
H_L(\al \cdot \al')=
H_L(\al)\circ H_L(\al'),
\end{equation}
where the dot denotes concatenation of $A$-paths.
\end{proposition}

\begin{proof}
Let $\al(t)$ be an $A$-path in $L$. We can find a
time-dependent section $\al_t$ of $A$ over $L$ such that
$\al_t(\gamma(t))=\al(t)$.
Using the notation above, we define a time-dependent section
$\tilde{\al}_t$ over the tubular neighborhood such that for $x\in L$
and $u\in p^{-1}(x)$
\[\tilde{\al}_t=\tilde{\al}_t^\parallel+\tilde{\al}_t^{\perp},
\quad \text{where } \tilde{\al}_t^\parallel\in E_{u},
\ \tilde{\al}_t^{\perp}\in (A_{F_{x}})_{u}.\]
The lifts $\tilde{\gamma}$ are the integral curves of
the vector field $X_t$ defined by
\[ X_t=\#\tilde{\al_t},\]
so $H_L(\al)_0$ is the map induced by the time-1 flow of $X_t$ on $F_{x_0}$.

The flow of $X_t$ is induced by the 1-parameter family of Lie
algebroid homomorphisms $\Phi_t^{\al_t}$ of $A$ obtained by
integrating the family $\tilde{\al}_t$ (see Section
\ref{sec:automorphisms}). The homomorphisms $\Phi_1^{\al_t}$ gives
a Lie algebroid isomorphism 
\[ H_L(\al):A_{F_{\gamma(0)}}\to A_{F_{\gamma(1)}},\] 
which covers $H_L(\al)_0$. Relation
(\ref{eq:holonomy:homomorphism}) follows since we have just shown
that $H_L(\al)$ is the time-1 map of some flow.
\end{proof}

We call $H_L(\al)$ the \emph{$A$-holonomy} of the
$A$-path $\al(t)$. One extends the definition of $H_L$ for
piecewise smooth $A$-paths in the obvious way.

Denote by $\Agerm(A_{F_x})$ the group of germs at $0$ of Lie
algebroid automorphisms of $A_{F_x}$ which map $0$ to $0$, and
by $\Omega_A(L,x_0)$ the group of piecewise smooth
\emph{$A$-loops} based at $x_0$.

\begin{definition}
The \textsc{$A$-holonomy} of the leaf $L$ with base
point $x_0$ is the map
\[H_L:\Omega_A(L,x_0)\to\Agerm(A_{F_{x_0}}).\]
\end{definition}

Notice that the holonomy of a leaf $L$ depends on the tubular
neighborhood $\tilde{i}:\nu(L)\to M$, on the choice of
trivialization, and on the choice of complementary bundles.
However, two different choices lead to conjugate homomorphisms.

\begin{example} Suppose $\#$ is injective, so we have a regular foliation $\F$
and $A$ can be identified with $T\F$. Then $A_{F_x}$ is a trivial
Lie algebroid so $\Agerm(A_{F_{x_0}})$ can be identified with
$\Agerm(F_{x_0})$, the group of germs of diffeomorphisms of
$F_{x_0}$ which map 0 to 0. Also, $E= A$ so the horizontal lift
$h(u,\al)$ is the unique tangent vector to the leaf through $u$
which projects to $\#\al$. We conclude that for a regular
foliation, the Lie algebroid holonomy coincides with the usual
holonomy.
\end{example}

\begin{example}
Let $A=T^*M$ be the Lie algebroid of a Poisson manifold $(M,\Pi)$.
In this case there is a natural choice for the complementary
subbundle $E$, namely $E_u=(T_u F_x)^0$. It follows from the
results in \cite{Fernandes:article:1} that, for this choice, each
automorphism $H_L(\al)$ covers a Poisson automorphism, and
is in fact determined by the Poisson automorphism that it covers.
Therefore, in this case, the Lie algebroid holonomy homomorphism is
essentially the same as the Poisson holonomy defined in
\cite{Fernandes:article:1}.
\end{example}

\begin{example}
\label{ex:holonomy:not:invariant}
Let $A=M\times\gg$ be the transformation Lie algebroid associated with
some infinitesimal action $\rho:\gg\to\X^1(M)$. Given $\al\in \gg$
we can identify it with a constant section of $A$. Let $x_0\in M$ be a
fixed point for the action and take $L=\set{x_0}$. Then $F_{x_0}=M$
and $A_{F_0}=A$, so $E$ is the trivial bundle over $M$. The horizontal
lift is given by $h(u,\al)=\#_u\al=\rho(\al)\cdot u$. If $\Psi:G\times M\to M$
denotes some local Lie group action that integrates $\rho$, we find
for the constant $A$-path $\al(t)=\al$
\[ H_L(\al)(x,v)=(\Psi(\exp(\al),x),Ad(\exp(\al))\cdot v),\]
which is a Lie algebroid automorphism of $A_{F_{x_0}}\simeq
M\times\gg$. Note that an $A$-loop with base path homotopic to a constant path 
might have non-trivial holonomy.
\end{example}

\subsection{Reduced Holonomy}
The Lie algebroid holonomy defined in the previous sections is not
homotopy invariant (Example \ref{ex:holonomy:not:invariant}).
Following the constructions given in \cite{Fernandes:article:1} and
\cite{Ginzburg:article:1} for the Poisson case, we can give a notion
of \emph{reduced holonomy} which is homotopy invariant.

Recall that $\Agerm(A_{F_x})$ denotes the group of germs at $0$ of Lie
algebroid automorphisms of $A_{F_x}$ which map $0$ to $0$.  We shall
denote by $\Ogerm(A_{F_x})$ the corresponding group of germs of outer Lie
algebroid automorphisms (see the end of Section \ref{sec:automorphisms}).

\begin{proposition}
\label{prop:reduced:holonomy}
Let $x\in L\subset M$ be a leaf of $A$ with associated $A$-holonomy
$H_L:\Omega_A(L,x)\to\Agerm(F_{x})$. If $\al_1(t)$ and
$\al_2(t)$ are $A$-loops based at $x$ with base paths
$\gamma_1\sim\gamma_2$ homotopic then $H_L(\al_1)$ and
$H_L(\al_2)$ represent the same equivalence class in $\Ogerm(F_{x})$.
\end{proposition}

\begin{proof}
Recall that any piecewise smooth path $\gamma\subset L$ can be made into
an A-path. By Proposition \ref{prop:holonomy}, property
(\ref{eq:holonomy:homomorphism}),
it is enough to show that for every
$x\in L$ there exists a neighborhood $U$ of $x$ in $L$ such that if
$\gamma(t)\subset U$ is a piecewise smooth loop based at $x$ and
$\al(t)\in A$ is a piecewise smooth family with
$\#\al=\dot{\gamma}$ then $H_L(\al)$ is a inner automorphism of
$A_{F_x}$.

We use the same notation as in the proof of Proposition
\ref{prop:holonomy}, so we construct a time-dependent
section $\tilde{\al}_t$ in a tubular neighborhood of $L$ which
decomposes as $\tilde{\al}_t=\tilde{\al}_t^\parallel+
\tilde{\al}_t^{\perp}$, and $H_L(\al)$ is obtained by
integrating this section up to time 1.

It is clear that the parallel component
$\#\tilde{\al}^{\parallel}_t$ has no effect on
the holonomy. Hence we can assume that $L=\set{x}$, $F_x=M$,
$\gamma$ is a constant path and $\tilde{\al}_t=\tilde{\al}^{\perp}_t$.
But then $\Phi_t^{\tilde{\al}_t}$ is a 1-parameter family of
automorphisms of $A_{F_x}$ with
$\Phi_1^{\tilde{\al}_t}=H_L(\al)$, so
we conclude that $H_L(\al)$ is an inner automorphism of $A_{F_x}$.
\end{proof}

Given a loop $\gamma$ in a leaf $L$ we shall denote by
$\bar{H}_L(\gamma)\in\Ogerm(A_{F_x})$ the equivalence class of
$H_L(\al)$ for some piece-wise smooth family $\al(t)$
with $\#\al(t)=\gamma(t)$. The map $\bar{H}_L:\Omega(L,x)\to\Ogerm(A_{F_x})$
will be called the \emph{reduced holonomy homomorphism} of $L$.
This maps extends to continuous loops and, by a standard argument, it
induces a homomorphism $\bar{H}_L:\pi_1(L,x)\to\Ogerm(A_{F_x})$ where
$\pi_1(L,x)$ is the fundamental group of $L$ (the use of the same letter
to denote both these maps should not cause any confusion).

\subsection{Stability}
Recall that, for a foliation $\F$ of a manifold $M$, a
\emph{saturated set} is a set $S\subset M$ which is a union
of leaves of $\F$. A leaf $L$ is called
\emph{stable} if it has arbitrarily small saturated neighborhoods.
In the case of the orbit foliation of a Lie algebroid a set is
saturated iff it is invariant under all inner automorphisms. Hence,
a leaf is stable iff it is has arbitrarily small neighborhoods
which are invariant under all inner automorphisms.

We shall call a leaf $L$ \emph{transversally stable} if $N\cap L$
is a stable leaf for the transverse Lie algebroid structure $A_N$,
i.e., if $N$ has arbitrarily small neighborhoods of $N\cap L$
which are invariant under all inner automorphisms of $A_N$.

The following result is a generalization of the Reeb Stability Theorem for
regular foliations.
\begin{theorem}[Stability Theorem]
\label{thm:local:stability}
Let $L$ be a compact, transversally stable leaf,
with finite reduced holonomy. Then $L$ is stable, i.e., $L$ has
arbitrarily small neighborhoods which are invariant under all
inner automorphisms. Moreover, each leaf near $L$ is a bundle over $L$
with fiber a finite union of leaves of the transverse Lie
algebroid structure.
\end{theorem}

\begin{proof}
Assume first that $L$ has trivial reduced holonomy and fix a base
point $x_0\in L$. We choose an embedding of $p:\nu(L)\to L$ in $M$, a
complementary subbundle $E$ and trivialization so we can define the
holonomy map $H_L$. Also, we choose a Riemannian metric on
$L$.

By compactness of $L$, there exists a number $c>0$ such that every
point $x\in L$ can be connected to $x_0$ by a smooth $A$-path
of length $<c$. For some inner product on $F_{x_0}$, let
$D_\eps$ be the disk of radius $\eps$ centered at $0$. For each
$\eps>0$, we can choose a neighborhood $U\subset D_\eps$ such that:
\begin{enumerate}
\item[(i)] for any piecewise-smooth $A$-path in $L$,
starting at $x_0$, with length $\le 2c$ and for any $u\in U$,
there exists a lifting with initial point $u$;
\item[(ii)] the lifting of any $A$-loop based at $x_0$
with initial point $u\in U$ has end point in $U$;
\item[(iii)] $U$ is invariant under all inner automorphisms of $A_{F_{x_0}}$;
\end{enumerate}
In fact, let $\al_1,\dots,\al_k$ be $A$-loops such that their base
loops $\gamma_1,\dots,\gamma_k$ are generators of
$\pi_1(L,x_0)$, and let $\Phi_i$ be Lie algebroid automorphisms
which represent the germs $H_L(\al_i)$. Since the reduced
holonomy is trivial, there is a neighborhood $U'$ of $0$ in
$F_{x_0}$ such that
$U\subset\text{domain}(\phi_1)\cap\cdots\cap\text{domain}(\phi_k)$,
and $\Phi_i|U'\in\Inn(A_{F_{x_0}})$, for all i. Since $L$ is
transversally stable, we can choose a smaller neighborhood $U\subset
U'$ invariant under all inner automorphisms.

Given $x\in L$ and an $A$-path $\al(t)$ connecting $x_0$ to $x$,
let us denote by $\sigma_{\al}:U\to F_{x}$ the diffeomorphism
defined by lifting. It follows from i) and ii) above that if
$\al'(t)$ is an $A$-path homotopic to $\al(t)$ then
$\sigma_{\al}(U)=\sigma_{\al'}(U)$. It follows from
iii) that $\sigma_{\al}(U)$ is also invariant under all inner
automorphisms.

Let $V$ be a neighborhood of $L$ in $M$. There exists $\eps(x)>0$
such that for the corresponding $U_x\subset D_{\eps(x)}$ we have
$\sigma_{\al}(U_x)\subset V\cap F_x$. By compactness of
$L$, we can choose $\eps>0$ (independent of $x\in L$) such that for
the corresponding $U\subset D_{\eps}$ we have
\[ \sigma_{\al}(U)\subset V\cap F_x\]
Set
\[ V_0=\bigcup_{\al}\sigma_{\al}(U).\]
Then $V_0\subset V$ is a open neighborhood of $L$ which is
invariant under all inner automorphisms of $M$. Therefore, $L$ is stable.

If $u,u'\in V_0$ are two points in the same leaf of $A$ such that
$p(u)=p(u')=x$, then there is a path $\tilde{\gamma}$ in this
leaf connecting these two points. We can choose a loop $\al(t)$
in $L$ based at $x$ such that $\tilde{\gamma}$ is a horizontal lift
of this loop.
Thus $u'$ is the image of $u$ by $H_L(\al)$ which is a inner
automorphism of $V_0\cap F_x$. Therefore, $u$ and $u'$ lie in the same
leaf of $V_0\cap F_x$. We conclude that each leaf of $M$ near $L$
is a bundle over $L$ with fiber a leaf of the transverse
Lie algebroid structure.

Assume now that $L$ has finite reduced holonomy. We let
$q:\tilde{L}\to L$ be a finite covering space such that
$q_*\pi_1(\tilde{L})=\Ker\bar{H}_L\subset\pi_1(L)$. If we embed
$\nu(L)$ into $M$ as above, and let $\nu(\tilde{L})$ be the pull
back bundle of $\nu(L)$ over $\tilde{L}$, we have a unique Lie
algebroid structure $\tilde{A}$ over $\nu(\tilde{L})$ and a Lie
algebroid homomorphism $\Phi:A \to\tilde{A}$ which covers the
natural map $\nu(\tilde{L})\to\nu(L)$. Moreover, the reduced
holonomy of $\tilde{A}\to\nu(\tilde{L})$ along $\tilde{L}$ is
trivial, so we can apply the above argument to $\nu(\tilde{L})$ and
the theorem follows.
\end{proof}

\begin{demo}{Remark}
If a leaf $L$ is transversally stable and $x\in L$, let $N$ denote a
stable neighborhood of $F_x$. For each $A$-loop $\al$, the
holonomy $H_L(\al)$ induces a homeomorphism of the
orbit space of $N$, for the transverse Lie algebroid structure,
mapping zero to zero. If $\al_1(t)$ and $\al_2(t)$
are $A$-loops such that $H_L(\al_1)$ and
$H_L(\al_2)$ represent the same class in $\Ogerm(F_x)$,
then they induce the same germ of homeomorphism of the orbit space
mapping zero to zero. In \cite{Dazord:article:1} holonomy of a
general, transversally stable, foliation is defined using germs of
homeomorphisms of the orbit space, which in the case of a foliation
defined by a Lie algebroid coincide with these homeomorphisms.
\end{demo}

\subsection{Linear Holonomy}
\label{sec:linear:holonomy}
Let $\pi:A\to M$ be a Lie algebroid and $i:L\hookrightarrow M$ a
leaf of $M$ with holonomy
$H_L:\Omega_A(L,x)\to\Agerm(F_x)$ (once appropriate data has been
fixed). Over $T_0 F_x\simeq F_x$ we consider the Lie algebroid
$A^\L_{F_x}\simeq \gg\times F_x$, where $\gg=\Ker\#_x$ which is the
linear approximation at $0$ to the transverse Lie algebroid
structure $A_{F_x}$. Also, we denote by $\Aut(A^\L_{F_x})$ the set
of linear Lie algebroid automorphisms of $A^\L_{F_x}$. There is a
map $d:\Agerm(A_{F_x})\to\Aut(A^\L_{F_x})$ which assigns to a germ
of a Lie algebroid automorphism of $A_{F_x}$, mapping zero to zero,
its linear approximation at 0. Obviously, we can identify such a
linear map with a pair $(\phi,\psi)$ where $\phi\in \textrm{GL}(F_x)$ and
$\psi$ is a Lie algebra automorphism of $\Ker\#_x=\gg$.

\begin{definition}
The \textsc{linear $A$-holonomy} of the leaf $L$ with base
point $x_0$ is the map
\[H_L^\L\equiv dH_L:\Omega_A(L,x)\to \Aut(\gg)\times \textrm{GL}(F_x).\]
\end{definition}

One can also define the \emph{reduced linear $A$-holonomy
homomorphism} of a leaf $L$ to be the class of $H_L^\L(\gamma,\al)$ in
$\Out(\gg)\times \textrm{GL}(F_x)$ where
$\Out(\gg)=\Aut(\gg)/\Inn(\gg)$. The reduced holonomy is homotopy
invariant. For Poisson manifolds, linear Poisson holonomy was first
introduced by Ginzburg and Golubev in \cite{Ginzburg:article:1}.

There is an alternative approach to linear holonomy using a linear
$A$-connection which generalizes the Bott connection of ordinary
foliation theory. In this differential operator formulation linear
holonomy arises as the holonomy of a Lie algebroid connection.

We consider first the vector bundle $\Ker\#|_L$ over the leaf $L$
where we have a \emph{Bott $A$-connection} defined as follows:
Given an element $\al\in A_x$, where $x\in L$, and a section
$\beta$ of $\Ker\#|_L$, we take (local) sections
$\tilde{\al},\tilde{\beta}\in\Gamma(A)$ such that
$\tilde{\al}_x=\al$, $\tilde{\beta}|_L=\beta$, and we set:
\begin{equation}
\label{Bott:connection}
\nabla^L_\al\beta\equiv {[\tilde{\al},\tilde{\beta}]}|_x.
\end{equation}

\begin{lemma}
$\nabla^L$ associates to each section $\al$ of $A$ along $L$ a
linear operator $\nabla_\al:\Gamma(\Ker\#|_L)\to\Gamma(\Ker\#|_L)$.
\end{lemma}
\begin{proof}
To check that expression (\ref{Bott:connection}) is independent of the
extensions considered, we fix a local basis of sections
$\set{\gamma^1,\dots,\gamma^r}$ for $A$ in a neighborhood of $x$. If
we write
\[\tilde{\al}=\sum_s \tilde{a}_s \gamma^s, \quad
\tilde{\beta}=\sum_t \tilde{b}_t \gamma^t,\]
for some functions $\tilde{a}_s$ and $\tilde{b}_t$, we
compute
\[ \nabla^L_\al\beta=\sum_{s,t}\left(\tilde{a}_s(x)\tilde{b}_t(x)
[\gamma^s,\gamma^t]|_x +
\tilde{a}_s(x)\#\gamma^s(\tilde{b}_t)|_x\gamma^t|_x\right).\]
This expression shows that $\nabla^L_\al\beta$ only depends on the value
of $\tilde{\al}$ at $x$ and the values of $\tilde{\beta}$ along
$L$, i.e., $\al$ and $\beta$.

Relation (\ref{Bott:connection}) also shows that $\nabla^L_\al \beta$ is in
the kernel of $\#$ and so is a section of $\Ker\#|_L$.
\end{proof}

Next we consider the conormal bundle $\nu^*(L)=\set{\omega\in
T^*_LM: \omega|TL=0}$ over the leaf $L$ where we also have a
\emph{Bott $A$-connection} defined as follows: Given an element
$\al\in A_x$, where $x\in L$, and a section $\omega$ of $\nu^*(L)$,
we take a section $\tilde{\al}\in\Gamma(A)$ and a 1-form
$\tilde{\omega}\in \Omega^1(M)$ such that $\tilde{\al}_x=\al$,
$\tilde{\omega}|_L=\omega$, and we set:
\begin{equation}
\label{Bott:connection:*}
\check{\nabla}^L_\al\omega\equiv \Lie_{\#\tilde{\al}}\tilde{\omega}|_x.
\end{equation}

A proof similar to the lemma above shows that:
\begin{lemma}
$\check{\nabla}^L$ associates to each section $\al$ of $A$ along $L$ a
linear operator
$\check{\nabla}_\al:\Gamma(\nu^*(L))\to\Gamma(\nu^*(L))$.
\end{lemma}

It is also easy to check that $\nabla^L$ and $\check{\nabla}^L$ satisfy the
analogue of properties (i) to (iv) of Proposition
\ref{prop:derivative:properties:2}. Note however that, in general,
$\nabla^L$ and $\check{\nabla}^L$ \emph{do not give} genuine Lie algebroid
connections since they are only defined for sections of $A$ along
$L$.

It is convenient to consider the connections $\nabla^L$ and
$\check{\nabla}^L$ all together, rather than leaf by leaf, so we set:

\begin{definition}
\label{defn:basic:connection}
A linear connection $\nabla$ on $A$ is called a
\textsc{basic connection} if
\begin{enumerate}
\item[(i)] $\nabla$ is compatible with the Lie algebroid structure, i.e.,
there exists a linear $A$-connection $\check{\nabla}$ on $TM$ such that
\[\check{\nabla}\#=\#\nabla.\]
\item[(ii)] $\nabla$ restricts to $\nabla^L$ on each leaf $L$, i.e.,
if $\al,\beta\in\Gamma(A)$ with $\#\beta|_L=0$ then
\[ \nabla_\al\beta |_L=\nabla^L_\al\beta.\]
\item[(iii)] $\check{\nabla}$ restricts to $\check{\nabla}^L$ on each
leaf $L$, i.e., if $\al\in\Gamma(A)$ and $\omega\in\Omega^1(M)$ with
$\omega|_{TL}=0$ then
\[ \nabla_\al\omega |_L=\check{\nabla}^L_\al\omega.\]
\end{enumerate}
\end{definition}

The holonomy along a leaf $L$ of a basic connection $\nabla$ coincides
with the linear holonomy of $L$ introduced above: the holonomy of
the basic connection $\nabla$ determines endomorphisms of the fiber
$A_x$ which map $\ker\#_x$ isomorphically into itself, and these
are the linear holonomy maps.

\begin{proposition}
Every Lie algebroid has a basic connection.
\end{proposition}

\begin{proof}
In fact, let us see that the compatible connection $\nabla$
constructed in the proof
of Proposition \ref{prop:compatible:connection} is a basic connection.
We use the same notation as in that proof, so if $L$ is a leaf of $M$ and
$\#\beta|_L=0$, we write $\beta=\sum_t b_t\al^t$ and we have
\begin{align*}
\nabla_{\al^s}^{(a)}\beta |_L&=\sum_t \left(b_t \nabla_{\al^s}^{(a)}\al^t+
\#\al^s(b_t)\al^t\right)\\
                      &=\sum_t \left(\sum_u b_t c^{st}_u\al^u+
\sum_j b^{sj}\frac{\partial b_t}{\partial x^j}\al^t\right)
=[\al^s,\beta]|_L
\end{align*}
Therefore, for any 1-form $\al=\sum_s  a_s \al^s$, we get
\begin{align*}
\nabla^{(a)}_\al\beta |_L&=\sum_s  a_s \nabla^{(a)}_{\al^s}\beta |_L\\
&= \sum_s  a_s [\al^s,\beta]|_L=[\al,\beta]|_L,
\end{align*}
since $\#\beta|_L=0$. It follows that for any 1-form $\al$ we have
\[ \nabla_\al\beta |_L=\sum_a \phi_a \nabla^{(a)}_\al\beta |_L=\nabla^L_\al\beta.\]

Similarly, for the connection $\check{\nabla}$, we have
\[ \check{\nabla}^{(a)}_{\al^s}dx^i=
\sum_j \frac{\partial b^{si}}{\partial x^j}dx^j=\Lie_{\#\al^s}dx^i,\]
so if $\omega|_{TL}=0$ we find
$\check{\nabla}^{(a)}_{\al}\omega|_L=\Lie_{\#\al}\omega$, and it follows
that
\[ \check{\nabla}_{\al}\omega|_L=\sum_a \phi_a \check{\nabla}^{(a)}_{\al}\omega|_L
=\check{\nabla}^{L}_{\al}\omega.\]
Since $\nabla\#=\#\check{\nabla}$ we conclude that $\nabla$ defines a basic
connection.
\end{proof}

In the theory of regular foliations basic connections arise as
connections in the normal bundle of the foliation (or equivalently
in the conormal bundle). In the case of a Lie algebroid, as was
first argued in \cite{Evens:article:1}, the ``formal difference''
$E=A\ominus T^*M$ plays the role of the (co)normal bundle. Now,
a basic connection $\nabla$ in $A$ gives an $A$-connection
$\tilde{\nabla}$ in the bundle $E$ (see Section
\ref{sec:compatible}):
\[ \tilde{\nabla}_\al(\beta+\omega)=
\nabla_\al\beta+ \check{\nabla}_\al\omega, \qquad
\al,\beta\in\Gamma(A),\ \omega\in\Omega^1(M).\]
We shall see later that this connection plays a key role in
defining
\emph{secondary characteristic classes} for Lie algebroids.

The following result shows that the curvature $\tilde{R}$ of
$\tilde{\nabla}$ vanishes along a leaf.
\begin{proposition}
\label{prop:basic:connection}
Let $\nabla$ be a basic connection and $L$ a leaf of $A$. Denote by $R$ and
$\check{R}$ the curvature of the connections $\nabla$ and $\check{\nabla}$.
If $\gamma$ is a section of $A$ such that $\#\gamma|_L=0$, then
\[ R(\al,\beta)\gamma|_L=0.\]
Similarly, if $\omega\in\Omega^1(M)$ is a differential form such that
$\omega|_{TL}=0$ then
\[ \check{R}(\al,\beta)\omega|_L=0.\]
\end{proposition}

\begin{proof}
If $\nabla$ is any basic connection and $\#\gamma|_L=0$, we have
$\nabla_\al\gamma |_L=[\al,\gamma]|_L$, so
expression (\ref{eq:curvature:derivative}) for the curvature tensor, gives
\[ R(\al,\beta)\gamma|_L=[\al,[\beta,\gamma]]|_L-[\beta,[\al,\gamma]]|_L-
[[\al,\beta],\gamma]|_L.\]
But the right hand side is zero, because of the Jacobi identity.

Similarly, if $\omega\in\Omega^1(M)$ is a differential form such that
$\omega|_{TL}=0$, we have
$\check{\nabla}_\al\omega|_L=\Lie_{\#\al}\gamma|_L$.
Hence, using $\#[\al,\beta]=[\#\al,\#\beta]$ and the well known formula
for the Lie derivative of the Lie bracket of vector fields, we  find
\[ R(\al,\beta)\omega|_L=\Lie_{\#\al}(\Lie_{\#\beta}\gamma)|_L-
\Lie_{\#\beta}(\Lie_{\#\al}\gamma)|_L-
\Lie_{[\#\al,\#\beta]}\gamma|_L=0,\]
so the second relation also holds.
\end{proof}

\begin{demo}{Remark}
Although the curvature of a basic connection vanishes along $\ker\#$,
the holonomy along $\#$ need not be discrete (this is because of the
presence of an extra term in the Holonomy Theorem
\ref{thm:holonomy}). Hence, in general, linear holonomy is not
discrete and also not homotopy invariant (cf.~Example
\ref{ex:holonomy:not:invariant}). However, if one can find a basic
\emph{$\F$-connection} then one gets discrete holonomy.  Such is the
case whenever $\#$ is injective, so the orbit foliation is
regular, and (linear) holonomy coincides with
standard (linear) holonomy of a regular foliation.
\end{demo}

\section{Characteristic Classes}

\subsection{Chern-Weil Homomorphism}

The usual Chern-Weil theory for characteristic classes extends to
$A$-connections. For contravariant connections this was already
discussed in \cite{Fernandes:article:1,Vaisman:art:1}. For general Lie
algebroids the theory is similar and only a short account will be
given as we shall need it later in the section when
we discuss secondary characteristic classes.

We consider a principal G-bundle $p:P\to M$ furnished with an
$A$-connection. Given any symmetric, $\Ad(G)$-invariant,
$k$-multilinear function
\[ P:\gg\times\cdots\times\gg\to\Rr\]
we can define a $2k$-section $\lambda(P)$ of $A$ as
follows. If $U_j$ is a trivializing neighborhood, and
$\al_1,\dots,\al_{2k}$ are sections of $A$ over $U_j$ then we set
\begin{multline}
\label{eq:Chern-Weil:homomorphism}
\lambda(P)(\al_1,\dots,\al_{2k})=\\
\sum_{\sigma\in S_{2k}}
(-1)^{\sigma}
P(\Omega_j(\al_{\sigma(1)},\al_{\sigma(2)}),\dots,\Omega_j(\al_{\sigma(2k-1)},\al_{\sigma(2k)})),
\end{multline}
where $\set{\Omega_j}$ are local curvature 2-sections. By the
transformation rule for the local curvature 2-sections, this formula
actually defines a $2k$-section $\lambda(P)\in\Gamma(\wedge^{2k}A^*)$ on the
whole of $M$.

\begin{proposition}
For any symmetric, invariant, $k$-multilinear function $P$, the
$2k$-section $\lambda(P)$ is closed:
\begin{equation}
d_A \lambda(P)=0.
\end{equation}
\end{proposition}

\begin{proof}
We compute
\begin{align*}
d_A \lambda(P)&=
k P(d_A \Omega_j,\dots,\Omega_j)\\
&=k P(d_A \Omega_j+[\omega_j,\Omega_j],\dots,\Omega_j)=0,
\end{align*}
where we have used first the linearity and symmetry of $P$, then the
$\Ad(G)$-invariance of $P$, and last the Bianchi identity.
\end{proof}

Therefore, to an invariant, symmetric, $k$-multilinear function
$P\in I^k(G)$ we can associate an $A$-cohomology class
$[\lambda(P)]\in H^{2k}(A)$, and in fact we have:

\begin{proposition}
\label{prop:connection:invariance:1}
The cohomology class $[\lambda(P)]$ is independent of the
$A$-connec\-tion used to define it.
\end{proposition}

\begin{proof}
Suppose we have two $A$-connections in $P(M,G)$ with
connection 1-sections $\omega^0$ and $\omega^1$,
and denote by $\lambda^0(P)$ and $\lambda^1(P)$ the $2k$-sections they
define through (\ref{eq:Chern-Weil:homomorphism}).  We construct a
1-parameter family of connections with connection 1-section
$\omega^t=t\omega^1+(1-t)\omega^0$, $t\in [0,1]$, and we
denote by $\Omega^t$ its curvature 2-section.

By the transformation rule (\ref{eq:transform:connection}) for the
local connection 1-sections, the difference
$\omega^{1,0}_j=\omega^1_j-\omega^0_j$ is a $\gg$-valued
1-section, and we get a well defined $(2k-1)$-section
$\lambda^{1,0}(P)$ by setting
\begin{multline}
\label{eq:invariants:2}
\lambda^{1,0}(P)(\al_1,\dots,\al_{2k-1})=\\
\sum_{\sigma} C_\sigma \int_{0}^1
P(\omega^{1,0}_j(\al_{\sigma(1)}),\Omega^t_j(\al_{\sigma(2)},\al_{\sigma(3)}),
\dots,\Omega^t_j(\al_{\sigma(2k-2)},\al_{\sigma(2k-1)}))dt.
\end{multline}
where $C_\sigma=k (-1)^{\sigma}$, and the sum is over all permutations 
in $S_{2k-1}$.
We claim that
\begin{equation}
\label{eq:delta:lambda:2}
d_A\lambda^{1,0}(P)=\lambda^1(P)-\lambda^0(P),
\end{equation}
so $[\lambda^1(P)]=[\lambda^0(P)]$.

To prove (\ref{eq:delta:lambda:2}), we note that if we differentiate
the structure equation (\ref{eq:first:structure:equation}) we obtain
\begin{equation}
\label{eq:tderivative:curvature}
\frac{d}{dt}\Omega^t_j=d_A  \omega^{1,0}_j+[\omega^t_j,\omega^{1,0}_j].
\end{equation}
Hence, using Bianchi's identity, we have
\begin{align*}
kd_A&\int_0^1 P(\omega^{1,0}_j,\Omega^t_j,\dots,\Omega^t_j)dt=\\
&=k\int_0^1 P(d_A \omega^{1,0}_j,\Omega^t_j,\dots,\Omega^t_j)\\
&\qquad \qquad \qquad
+P(\omega^{1,0}_j,d_A \Omega^t_j,\dots,\Omega^t_j)
+P(\omega^{1,0}_j,\Omega^t_j,\dots,d_A \Omega^t_j)dt\\
&=
k\int_0^1
P(\frac{d}{dt}\Omega^t_j-[\omega^t_j,\omega^{1,0}_j],\Omega^t_j,\dots,\Omega^t_j)\\
&\qquad \qquad \qquad
-P(\omega^{1,0}_j,[\omega^t_j,\Omega^t_j],\dots,\Omega^t_j)
-P(\omega^{1,0}_j,\Omega^t_j,\dots,[\omega^t_j,\Omega^t_j])dt\\
&=
k\int_0^1
P(\frac{d}{dt}\Omega^t_j,\Omega^t_j,\dots,\Omega^t_j)dt\\
&=
\int_0^1
\frac{d}{dt}P(\Omega^t_j,\Omega^t_j,\dots,\Omega^t_j)dt=
P(\Omega^1_j,\dots,\Omega^1_j)-P(\Omega^0_j,\dots,\Omega^0_j),
\end{align*}
so the claim follows.
\end{proof}

If we set
\[ I^\bullet(G)=\bigoplus_{k\ge 0}I^k(G),\]
the assignment $P\mapsto [\lambda(P)]$ gives a map
$I^\bullet(G)\to H^\bullet(A)$, which is in fact a ring
homomorphism, and which we call the \emph{$A$-Chern-Weil
homomorphism} of the Lie algebroid. The fact that this map is a
ring homomorphism follows, for example, from the following proposition:

\begin{proposition}
\label{prop:A:Chern:Weil}
The following diagram commutes
\[
\xymatrix{
I^\bullet(G)\ar[r]\ar[dr]& H^\bullet_{\text{de Rham}}(M) \ar[d]^{\#^*} \\
& H^\bullet(A)}
\]
where on the top row we have the usual Chern-Weil homomorphism.
\end{proposition}

\begin{proof}
Choose an $A$-connection in $P$ which is induced by some covariant
connection. Given $P\in I^k(G)$, this covariant connection gives
a closed $(2k)$-form $\tilde{\lambda}(P)$ defined by a
formula analogous to (\ref{eq:Chern-Weil:homomorphism}), and which induces
the usual Chern-Weil homomorphism
$I^\bullet(G)\to H^\bullet_{\text{de Rham}}(M)$.
We check easily that
\[ \#^*\tilde{\lambda}(P)=\lambda(P),\]
so the proposition follows.
\end{proof}

Recall that the ring of invariant polynomials $I^\bullet(GL_q(\Rr))$ is generated by elements
$P_k\in I^k(GL_q(\Rr))$ such that $P_k(X,\dots,X)=\sigma_k(X)$,
where $\set{\sigma_1,\dots,\sigma_q}$
are the \emph{elementary symmetric functions} defined by:
\[ \det(\mu I-\frac{1}{2\pi}X)=\mu^q+\sigma_1(X)\mu^{q-1}+\cdots+\sigma_q(X).\]

Now consider a real vector bundle $p_E:E\to M$, with $\rank E=q$,
and let $p:P\to M$ be the associated principal
bundle with structure group $GL_q(\Rr)$. Choosing a Lie algebroid connection
on $P$ one defines the  \emph{kth $A$-Pontrjagin class} of $E$ as
\[ p_k(E,A)=[\lambda(P_{2k})]\in H^{4k}(A).\]
As usual, one does not need to consider the classes for odd $k$ since we have
\[ [\lambda(P_{2k-1})]=0,\]
as can be seen by choosing a connection compatible with a riemannian metric.
It is clear from Proposition \ref{prop:A:Chern:Weil} that
\[ p_k(E,A)=\#^*p_k(E).\]
where $p_k(E)$ are the usual Pontrjagin classes of $E$.

To compute these invariants one uses the $A$-derivative operator $\nabla$
on $E$, associated with the $A$-connection, and proceeds
as follows. For 1-sections $\al,\beta\in \Gamma(A)$, the curvature tensor $R$
defines a linear map $R_{\al,\beta}=R(\al,\beta):E_x\to E_x$ which satisfies
$R_{\al,\beta}=-R_{\beta,\al}$, and so $(\al,\beta)\to R_{\al,\beta}$ can be
considered as a $\gl(E)$-valued 2-section. By fixing a basis of local
sections for $E$, we have $E_x\simeq\Rr^q$ so we have
$R_{\al,\beta}\in\gl_q(\Rr)$.
(this matrix representation of $R_{\al,\beta}$ is defined only up to a change
of basis in $\Rr^q$). Hence, if
\[ P:\gl_q(\Rr)\times\cdots\times\gl_q(\Rr)\to\Rr\]
is a symmetric, $k$-multilinear function, $\Ad(GL_q(\Rr))$-invariant, we
a have a $2k$-section $\lambda(R)(P)$ defined by
\begin{equation}
\label{eq:Chern-Weil:homomorphism:connection}
\lambda(R)(P)(\al_1,\dots,\al_{2k})=\sum_{\sigma\in S_{2k}} (-1)^{\sigma}
P(R_{\al_{\sigma(1),\sigma(2)}},\dots,R_{\al_{\sigma(2k-1),\sigma(2k)}}).
\end{equation}
It is easy to see that $\lambda(P)=\lambda(R)(P)$, so this gives a
procedure to compute the $A$-Chern-Weil homomorphism and the
$A$-Pontrjagin classes.

Similar considerations apply to other characteristic classes. One can
define, e.g., the $A$-Chern classes $c_k(E,A)$ of a complex
vector bundle $E$ and they are just the images by
$\#^*$ of the usual Chern classes of $E$.

The fact that all these classes arise as image by $\#^*$ of well known
classes is perhaps a bit disappointing. However, we shall see below
that one can define secondary characteristic classes which are true
invariants of the Lie algebroid, in the sense that they do not arise
as images by $\#^*$ of some de Rham cohomology classes. On the other hand,
the tangential Chern-Weil theory is usefull on its own, 
and have interesting applications already in the case of regular foliations 
(see e.g.~\cite{Moore:Schochet:book}).

\subsection{Secondary Characteristic Classes}

Whenever a form representing a (primary) characteristic class vanishes
one can introduce new (secondary) characteristic classes, a remark
that goes back to the original paper of Chern and Simons
(\cite{Chern:article:1}).  This remark was the starting point for the
theory of ``exotic'' characteristic classes for foliations (see
\cite{Bott:lectures:1}).

On a Lie algebroid a similar construction of ``exotic'' characteristic
classes can be done. This construction generalizes the
construction given in \cite{Fernandes:article:1} for the case of a
Poisson manifold, where it is was shown that Poisson secondary
characteristic classes give information on the topology, as well as,
the geometry of the symplectic foliation.

In the theory of (regular) foliations, the secondary characteristic
classes appear when we compare two connections on the normal bundle
each from a distinguished class. In the case of a Lie algebroid
the ``formal difference'' $E=A\ominus T^*M$ plays the role of the
(co)normal bundle, and again we compare two connections, each from
a distinguished class. So on a Lie algebroid $\pi:A\to M$, with
$\rank A=r$ and $\dim M=m$, we consider the following data:
\begin{enumerate}
\item[(i)] A basic connection with associated $A$-derivative operator
$\nabla$, so on $E$ we have an induced connection (see section
\ref{sec:linear:holonomy}) $\nabla^1=\nabla\oplus\check{\nabla}$;
\item[(ii)] A metric connection $\nabla^0$ on $E$, i.e., we have a covariant
connection $\bar{\nabla}^0$ on $E$, which preserves some metric $g$ on
$E$, and we take $\nabla^0_\al=\bar{\nabla}^0_{\#\al}$;
\end{enumerate}

Given an invariant, symmetric, $k$-multilinear function
$P\in I^k(GL(r+m,\Rr))$ we consider the $(2k-1)$-section
$\lambda^{1,0}(P)\in \Gamma(\wedge^{2k-1}A^*)$ given by
(\ref{eq:invariants:2}).

\begin{proposition}
\label{prop:invariants:2:closed}
If $k$ is odd, $\lambda^{1,0}(P)$ is a $d_A$-closed
$(2k-1)$-section.
\end{proposition}

\begin{proof}
According to (\ref{eq:delta:lambda:2}) we have
\[
d_A\lambda^{1,0}(P)=\lambda^1(P)-\lambda^0(P).
\]
and we claim that $\lambda^1(P)=\lambda^0(P)=0$ if $k$ is odd (these
are the vanishing primary classes that we mentioned to above).

The proof that $\lambda^0(P)=0$ is standard: we can choose an
orthonormal basis of sections for $E$ so that the curvature 2-sections
take there values in $\mathfrak{so}(r+m,\Rr)$. But if
$X\in \mathfrak{so}(r+m,\Rr)$, we have $P_k(X)=0$
for any elementary symmetric function, since $k$ is odd. Hence we obtain
$\lambda^0(P)=0$.

Consider now the connection $\nabla^1$. Given $x\in M$ we choose local
coordinates $(x^j,y^j)$ around $x$ and a basis of sections
$\set{\al^1,\dots,\al^r}$ as in the Local Splitting Theorem. Then
$\set{dx^i,dy^j,\al^s}$ form a basis for $E$, and for the canonical
skew-symmetric bilinear form $(~,~)$ given by (\ref{eq:bilinear:form})
the only non-vanishing pairs are:
\[ (\al^i,dx^i)=1=-(dx^i,\al^i), \qquad
(\al^s,dy^j)=b^{sj}=-(dy^j,\al^s).\]
where $b^{sj}=b^{sj}(y)$. Since $\nabla^1$ is induced by a basic
connection, it is compatible with the Lie algebroid structure so from
proposition \ref{prop:compatible} we conclude:
\begin{align*}
(\nabla^1_\gamma \al^i,\al^l)&=-(\al^i,\nabla^1_\gamma \al^l),\\
(\nabla^1_\gamma \al^i,dx^l)&=-(\al^i,\nabla^1_\gamma dx^l), \quad (1\le i,l\le q)\\
(\nabla^1_\gamma dx^i,dx^l)&=-(dx^i,\nabla^1_\gamma dx^l).
\end{align*}
On the other hand, from Proposition \ref{prop:basic:connection} we find
\[R^1(\al,\beta)dy^j|_x=R^1(\al,\beta)\al^s|_x=0, \quad (j,s>q).\]
It follows that $R^1(\al,\beta)_x$ is represented in the basis
$(dx^i,\al^i,dy^j,\al^s)$ by a matrix of the form:
\begin{equation}
\label{eq:matrix:form}
\left(\begin{array}{cc}
B & 0 \\
C & 0
\end{array}\right),
\end{equation}
with $B$ a $(2q\times 2q)$ symplectic matrix. Now, if
$A$ is any matrix of this form, it is
clear that $\det(\mu I-A)=\det(\mu I-\tilde{A})$, where $\tilde{A}$ is the
same as $A$ with $C=0$, i.e., $\tilde{A}$ is symplectic. But if
$\tilde{A}$ is symplectic, we have $P_k(A)=0$
for any elementary symmetric function, since $k$ is odd. Hence $\lambda^1(P)=0$.
\end{proof}

Next we want to check that the cohomology class of
$\lambda^{1,0}(P)$ is independent of the connections used to
define it.

Given 3 connections with local connection 1-sections
$\omega^0_j,\omega^1_j,\omega^2_j$ we consider a family of
connections with local connection 1-sections 
$\omega^{s,t}_j=(1-s-t)\omega^0_j+s\omega^1_j+t\omega^2_j$,
where $(s,t)$ vary in the standard 2-simplex $\Delta_2$. We
introduce a $(2k-2)$-vector field $\lambda^{2,1,0}(P)$ by a
formula analogous to (\ref{eq:Chern-Weil:homomorphism:connection})
and (\ref{eq:invariants:2}):
\begin{equation}
\label{eq:invariants:3}
\lambda^{2,1,0}(P)=
k \sum_{\sigma\in S_{2k-2}} (-1)^{\sigma}\int_{\Delta_2}
P(\omega^{1,0}_j,\omega^{2,0}_j,\Omega^{s,t}_j,\dots,\Omega^{s,t}_j)dtds,
\end{equation}
and just like in the proof of Proposition
\ref{prop:connection:invariance:1}, one shows that
\begin{equation}
\label{eq:delta:lambda:3}
\delta\lambda^{2,1,0}(P)=\lambda^{1,0}(P)-
\lambda^{2,0}(P)+\lambda^{2,1}(P).
\end{equation}

\begin{proposition}
\label{prop:independency:connection}
The cohomology class $[\lambda^{1,0}(P)]$ is independent
of the connections used to define it.
\end{proposition}

\begin{proof}
Let $\nabla^1$ and $\tilde{\nabla}^1$ (resp.~$\nabla^0$ and
$\tilde{\nabla}^0$) be basic connections (resp.~riemannian connections).
It follows from (\ref{eq:delta:lambda:3}) that
\begin{multline*}
\lambda(\nabla^1,\nabla^0)(P)-
\lambda(\tilde{\nabla}^1,\tilde{\nabla}^0)(P)=
\delta\lambda(\tilde{\nabla}^1,\nabla^0,\tilde{\nabla}^0)(P)+
\delta\lambda(\tilde{\nabla}^1,\nabla^1,\nabla^0)(P)\\
-\lambda(\tilde{\nabla}^1,\nabla^1)(P)
-\lambda(\nabla^0,\tilde{\nabla}^0)(P).
\end{multline*}
Hence, it is enough to show that the cohomology classes of
$\lambda(\tilde{\nabla}^1,\nabla^1)(P)$ and
$\lambda(\nabla^0,\tilde{\nabla}^0)(P)$ are trivial.

Consider first the basic connections $\tilde{\nabla}^1$ and $\nabla^1$.
The linear combination $\nabla^{1,t}=(1-t)\nabla^1+t\tilde{\nabla}^1$ is also
a basic connection. If $x\in M$, we fix splitting coordinates $(x^i,y^j)$
around $x$ and sections $\set{\al^1,\dots,\al^r}$ as in the proof of
proposition \ref{prop:invariants:2:closed}.
Then we see that, with respect to the basis
$\set{dx^i,\al^i,dy^j,\al^s}$, the matrix
representations of $\tilde{\nabla}^1_\al-\nabla^1_\al$ and $R^t(\al,\beta)$
are of the form (\ref{eq:matrix:form}). Hence, we conclude that if
$P\in I^k(GL(m,\Rr))$, with $k$ odd,
\[P(\tilde{\nabla}^1_{\al_1}-\nabla^1_{\al_1},R^t(\al_2,\al_3),\dots,
R^t(\al_{2k-2},\al_{2k-1}))=0.\]
Therefore, $\lambda(\tilde{\nabla}^1,\nabla^1)(P)=0$, whenever
$\tilde{\nabla}^1$ and $\nabla^1$ are basic connections.

Now consider the riemannian connections $\nabla^0$ and $\tilde{\nabla}^0$. The
linear combination $\nabla^{0,t}=(1-t)\tilde{\nabla}^0+t\nabla^0$ is also a
riemannian connection. All these connections are induced from
covariant riemannian connections $\bar{\nabla}^0$, $\tilde{\bar{\nabla}}^0$ and
$\bar{\nabla}^{0,t}$, and we can define a differential form
$\lambda(\bar{\nabla}^0,\tilde{\bar{\nabla}}^0)(P)$ of degree $(2k-1)$ by a
formula analogous to (\ref{eq:invariants:2}). Moreover, this form is
closed (because $k$ is odd), and
$\#^*\lambda(\bar{\nabla}^0,\tilde{\bar{\nabla}}^0)(P)=\lambda(\nabla^0,\tilde{\nabla}^0)(P)$.
It follows from the homotopy invariance of
$H^\bullet_{\text{de Rham}}(M)$, using a
suspension argument, as in the usual
theory of secondary characteristic classes of foliations (see
\cite{Bott:lectures:1}, page 29), that
\[ [\lambda(\bar{\nabla}^0,\tilde{\bar{\nabla}}^0)(P)]=[\lambda(\bar{\nabla}^0,\bar{\nabla}^0)(P)].\]
Hence, the cohomology class
$[\lambda(\bar{\nabla}^0,\tilde{\bar{\nabla}}^0)(P)]$
vanishes and so does the class $[\lambda(\nabla^0,\tilde{\nabla}^0)(P)]$.
\end{proof}

\begin{demo}{Remark}
We have used a riemannian connection of the special form
$\nabla^0_\al=\bar{\nabla}^0_{\#\al}$. On the other hand,
in general, a Lie algebroid does not admit a compatible $A$-connection
of the form $\tilde{\nabla}_{\#\al}$. Hence, the basic connections are
``genuine'' $A$-connections, i.e., not induced by any covariant
connection.
\end{demo}
\vskip 10pt

\begin{definition}
The \emph{secondary characteristic classes}
$\set{m_k(A)}$ of a Lie algebroid are the Lie algebroid cohomology
classes
\begin{equation}
m_k(A)=[\lambda^{1,0}(P_{k})]\in H^{2k-1}(A),\qquad
(k=1,3,\dots).
\end{equation}
where $P_k$ are the elementary symmetric polynomials.
\end{definition}

Note that, by the remark above, these secondary characteristic classes are
``genuine'' Lie algebroid cohomology classes, i.e., they do not lie
in the image of $\#^*:H^\bullet_{\text{de Rham}}(M)\to H^\bullet(A)$
(see also the examples below where one can have trivial de Rham
cohomology and non-zero secondary characteristic classes).
\vskip 15pt

\begin{demo}{Remark}
In general, one can only define the characteristic classes $m_k$ for $k$ odd.
Assume, however, that $A$ admits flat riemannian connections and flat basic
connections. Then the proofs of Propositions
\ref{prop:invariants:2:closed} and \ref{prop:independency:connection} can be
carried through, in the class of flat connections, for \emph{any} $k$.
Hence, in this case, one can define characteristic classes $m_k$
for \emph{any} $k$.
\end{demo}

\subsection{The Modular Class}
\label{section:modular:class}
The modular class of a Lie algebroid was introduced in
\cite{Weinstein:article:2}, and further discussed in
\cite{Evens:article:1}. Extensions to more general algebraic settings
were given in \cite{Huebschmann:article:1,Huebschmann:article:2,Xu:article:1}.

Let us start by recalling the construction given in
\cite{Evens:article:1}. Consider the line bundle
$Q_A=\wedge^r A\otimes\wedge^m T^* M$. It is
easy to check that on this line bundle we have a flat $A$-connection $\nabla$
defined by:
\begin{multline}
\label{eq:connection:Lie}
\nabla_\al(\al^1\wedge\cdots\wedge\al^r\otimes \mu)=
\sum_{j=1}^r
\al^1\wedge\cdots\wedge[\al,\al^j]\wedge\cdots\wedge\al^r\otimes\mu+\\
\al^1\wedge\cdots\wedge\al^r\otimes\Lie_{\#\al}\mu,
\end{multline}
whenever $\al,\al^1,\dots,\al^r\in\Gamma(A)$ and $\mu\in\Gamma(\wedge^m T^*M)$.

Now assume first that $Q_A$ is trivial. Then we have a global section
$s\in\Gamma(Q_A)$ so that
\[\nabla_\al s=\theta_s(\al)s, \qquad \forall \al\in\Gamma(A).\]
Since $\nabla$ is flat, we see that $\theta_s$ defines a section of
$\Gamma(A^*)$ which is closed: $d_A\theta_s=0$. If $s'$ is another
global section in $\Gamma(Q_A)$, we have $s'=a s$ for some
nonvanishing smooth function $a$ on $M$, and we find
\[ \theta_{s'}=\theta_s+d_A\log|a|.\]
Therefore, we have a well defined cohomology class
\[ \text{\rm mod}\,(A)\equiv[\theta_s]\in H^1(A)\]
which is independent of the section $s$. If the line bundle $Q_A$ is
not trivial one considers the square $L=Q_A\otimes Q_A$, which is
trivial, and defines
\[ \text{\rm mod}\,(A)=\frac{1}{2}[\theta_s],\]
for some global section $s\in\Gamma(L)$.

\begin{definition}
The class $\text{\rm mod}\,(A)$ is called the \textsc{modular class} of the
Lie algebroid $A$.
\end{definition}

As was argued in \cite{Evens:article:1} one can
think of global sections of $Q_A$ (or $Q_A\otimes Q_A$) as
``transverse measures'' to $A$. The modular class is trivial iff
there exists a transverse measure which is invariant under the
flows of $X_\al\equiv(X_{f_\al},\#\al)$, for every section
$\al\in\Gamma(A)$. Hence, the modular class is an obstruction lying
in the first Lie algebroid cohomology group $H^1(A)$ to the
existence of a transverse invariant measure to $A$.

\begin{theorem}
For any Lie algebroid $A$
\begin{equation}
\label{eq:invariant:modular:class}
m_1(A)=\frac{1}{2\pi}\text{\rm mod\,}(A).
\end{equation}
\end{theorem}

\begin{proof}
Choose a basic connection $\nabla^1$ and a riemannian connection $\nabla^0$
relative to some metric on $E=A\ominus T^*M$. We consider the
transverse measure $s$ to $A$ associated with this metric. We claim that
\begin{equation}
\label{eq:invariant:modular:field}
\lambda^{1,0}(\tr)=\theta_s,
\end{equation}
so (\ref{eq:invariant:modular:class}) follows.

Observe that it is enough to show that
(\ref{eq:invariant:modular:field}) holds on the regular points of $M$,
since the set of regular points is an open dense set and both sides
are smooth sections in $\Gamma(A^*)$. So assume that $x\in M$ is a
regular point where $\rank=q$, and pick coordinates $(x^i)$ around $x$
and a basis of sections $\set{\al^1,\dots,\al^r}$ as in the Local Splitting
Theorem. Then $s$ is given locally by:
\[s=(\det g)^{\frac{1}{2}}\al^1\wedge\cdots\wedge
\al^r\otimes dx^1\wedge\cdots\wedge dx^m,\]
where $g=(g^{ij}(x))$ is the matrix of inner products formed by
elements in $\set{\al^s,dx^i}$.

As in the proofs of the previous section, one computes the trace of
the operator $\nabla^1_\al$ relative to the basis
$\set{\al^s,dx^i}$ to be
\[ \tr \nabla^1_{\al^s}=\sum_{u>q}c^{su}_u+\sum_{j>q}\frac{\partial
b^{sj}}{\partial x^j}.\]
Also, since $\nabla^0$ is a metric connection, we find:
\begin{align*}
0=&\bar{\nabla}^0_{\#\al^s}s\\
 =&\#\al^s((\det g)^{\frac{1}{2}})\al^1\wedge\cdots\wedge
\al^r\otimes dx^1\wedge\cdots\wedge dx^m+\\
 &\qquad+(\det g)^{\frac{1}{2}}(\bar{\nabla}^0_{\#\al^s}
\al^1\wedge\cdots\wedge\al^r\otimes dx^1\wedge\cdots\wedge dx^m+\\
&\qquad\qquad\qquad+\cdots+\al^1\wedge\cdots\wedge
\al^r\otimes dx^1\wedge\cdots\wedge \bar{\nabla}^0_{\#\al^s}dx^m)\\
 =&\left(\#\al^s((\det g)^{\frac{1}{2}})
 +(\det g)^{\frac{1}{2}}\tr
\bar{\nabla}^0_{\#\al^s}\right)\al^1\wedge\cdots\wedge\al^r\otimes
dx^1\wedge\cdots\wedge dx^m.
\end{align*}
So we conclude that:
\begin{multline}
\label{eq:modular:aux:1}
\tr (\nabla^1_{\al^s}-\nabla^0_{\al^s})s=
\#\al^s((\det g)^{\frac{1}{2}})
\al^1\wedge\cdots\wedge\al^r\otimes dx^1\wedge\cdots\wedge dx^m\\
+(\sum_{u>q}c^{su}_u+
\sum_{j>q}\frac{\partial b^{sj}}{\partial x^j})s.
\end{multline}

On the other hand, a straight forward computation using
(\ref{eq:connection:Lie}) and the various relations in the
Local Splitting Theorem at a regular point, shows that
\begin{align}
\label{eq:modular:aux:2}
\nabla_{\al^s}s&=\#\al^s((\det g)^{\frac{1}{2}})
\al^1\wedge\cdots\wedge\al^r\otimes dx^1\wedge\cdots\wedge dx^m\nonumber\\
&\qquad+(\det g)^{\frac{1}{2}}
[\al^s,\al^1]\wedge\cdots\wedge\al^r\otimes dx^1\wedge\cdots\wedge
dx^m\nonumber\\
&\qquad\qquad+\cdots+(\det g)^{\frac{1}{2}}
\al^1\wedge\cdots\wedge\al^r\otimes dx^1\wedge\cdots\wedge
\Lie_{\al^s}dx^m\nonumber\\
&=\#\al^s((\det g)^{\frac{1}{2}})
\al^1\wedge\cdots\wedge\al^r\otimes dx^1\wedge\cdots\wedge dx^m+
\nonumber\\
&\qquad \qquad \qquad\qquad \qquad \qquad
+(\sum_{u>q}c^{su}_u+\sum_{j>q}\frac{\partial b^{sj}}{\partial x^j})s.
\end{align}

Comparing (\ref{eq:modular:aux:1}) and (\ref{eq:modular:aux:2})
gives
\[ \nabla_\al s=\tr (\nabla^1_{\al}-\nabla^0_{\al})s,\]
so relation (\ref{eq:invariant:modular:field}) holds and the theorem follows.
\end{proof}

\begin{demo}{Remark}
In a recent preprint \cite{Crainic:preprint:1}, M.~Crainic proposes a
different approach to secondary characteristic classes for vector
bundles admitting flat $A$-connections. Let $E$ be any vector
bundle over $M$ which admits a flat$A$-connection and is trivial as
a vector bundle. Then for each symmetric, $G$-invariant,
$2k-1$-multilinear function $P$ one can define a $2k-1$-section
$\nu(P)\in\Gamma(\wedge^{2k-1}A^*)$ by setting
\begin{equation}
\nu(P)(\al_1,\dots,\al_{2k})=
\sum_{\sigma\in S_{2k}}
(-1)^{\sigma}
P(\omega_j(\al_{\sigma(1)}),\dots,\omega_j(\al_{\sigma(2k-1)})),
\end{equation}
where $\al_1,\dots,\al_{2k}$ are global sections of $A$. Here
$\omega_j$ denotes a connection local 1-section which, since we are assuming
that $E$ is trivial, is actually globally defined. One can check that
$\nu(P)$ is closed so actually defines a cohomology class in
$H^{2k-1}(A)$.

If $E$ is not trivial as vector bundle, one needs a kind of C\v{e}ch
cohomology argument given in \cite{Crainic:preprint:1} to define
these classes. These classes generalize the classes $\theta_s$
constructed above for a line bundle admitting a flat
$A$-connection.

The ``normal bundle'' $E=A\ominus T^*M$ in general admits no flat
$A$-connection (this is clear already in the Poisson case, where
$A=T^*M$) and hence it is not obvious if one can use this approach
to define the secondary classes $m_k(A)$. Perhaps, as stated in
\cite{Crainic:preprint:1}, it is possible to extend this approach to
so called \emph{representations up to homotopy} (in our notation,
\emph{flat connections up to homotopy}) so one can use the ``adjoint
representation'' of $A$ (see \cite{Evens:article:1} for details).
\end{demo}

\subsection{Examples} We now consider some of the classes of Lie
algebroids that we have mentioned in section 1 and compute their secondary
characteristic classes.

\paragraph{Regular Foliations}
Let $\F$ be a regular foliation, and denote by $A=T\F\subset TM$ the
associated integrable subbundle.
Observe that a section $\al\in\Gamma(A)$ is just a vector field on
$M$ tangent to $\F$, so $A$ is a Lie algebroid with anchor the
inclusion $A\subset TM$ and $[~,~]$ the usual Lie bracket of vector
fields.

First choose some riemannian connection in $M$ determining a splitting
\[ T^*(M)=A^*\oplus \nu^*(\mathcal{\F}),\]
where $\nu^*(\mathcal{\F})$) is the conormal bundle to the
foliation. We have an $A$-riemannian connection $\check{\nabla}^0$ such that:
\[ \check{\nabla}^0_\al(\beta+\gamma)=\nabla_{\#\al}^{0,\parallel}\beta+
\nabla_{\#\al}^{0,\perp}\gamma, \]
where $\beta$ and $\gamma$, are sections of $A^*=T^*(\mathcal{\F})$ and
$\nu^*(\mathcal{\F})$, and $\nabla^{0,\parallel}$ and
$\nabla^{0,\perp}$, are covariant riemannian connections in these
bundles.
We choose on $E=A\ominus T^*M$ the $A$-riemannian connection
$\nabla^0_\al=\nabla^{0,\parallel}_{\#\al}\oplus\check{\nabla}^0_\al$.
Note that we are using the same notation for a connection on a vector
bundle and on its dual, but in fact, taking $\nabla^0$ on $E$ is essentially
equivalent to take $\nabla^{0,\parallel}$ on $\nu^*(\mathcal{\F})$.

Now we take as a basic connection (in the sense of definition
\ref{defn:basic:connection}) a connection $\nabla^1$ on
$E\simeq A\oplus A^*\oplus \nu^*(\mathcal{\F})$ of the form
$\nabla^1=\nabla^{1,\parallel}\oplus \nabla^{1,\parallel}\oplus
\nabla^{1,\perp}$ where $\nabla^{1,\parallel}$ is a connection in $A$
and $\nabla^{1,\perp}$ is just a basic connection in
$\nu^*(\mathcal{\F})$ in the usual
sense of foliation theory (see \cite{Bott:lectures:1}, p. 33).
A straightforward computation shows that
\[
\lambda(\nabla^1,\nabla^0)(P)(\al_1,\dots,\al_{2k-1})
=\lambda(\nabla^{1,\perp},\nabla^{0,\perp})(P)((\#\al_1,\dots,\#\al_{2k-1}).
\]

Recall that in foliation theory (see \cite{Bott:lectures:1}, p.~66)
the forms
\begin{align*}
c_k&=\lambda(\nabla^{1,\perp})(\tilde{P}_{k}),\qquad (1\le k\le q)\\
h_{2k-1}&=\lambda(\nabla^{1,\perp},\nabla^{0,\perp})(\tilde{P}_{2k-1}),
\qquad (1\le 2k-1\le q),
\end{align*}
satisfy
\begin{align}
\label{eq:Gelfand:Fuks}
dc_k&=0,\qquad (1\le k\le \text{corank}(\F))\\
dh_{2k-1}&=c_{2k-1},\qquad (1\le 2k-1\le \text{corank}(\F)).
\end{align}
and so they can be used to define a homomorphism of graded algebras
\[ H^*(WO_q)\to H^*(M),\]
where $H^*(WO_q)$ is the relative Gelfand-Fuks cohomology of formal
vector fields in $\Rr^q$. This homomorphism is independent of the
connections and its image are the exotic or secondary characteristic
classes of foliation theory.

Observe that the $(2k-1)$-forms
$h_{2k-1}=\lambda(\nabla^{1,\perp},\nabla^{0,\perp})(P_{k})$
are not closed in general, but are closed along the leaves, so
its image under $\#$ is a closed $(2k-1)$-section of $T^*\F^*$. Hence,
$\#h_{2k-1}$ defines a tangential cohomology class, and one has
\begin{equation}
m_{2k-1}(T\F)=[\#h_{2k-1}]
\end{equation}
but, in general, $m_{2k-1}$ is not in the image of
$\#^*:H^\bullet_{\text{de Rham}}(M)\to H^\bullet(A)$. A simple
consequence of this relationship is that, for a regular foliation,
the characteristic classes $m_k(T\F)$ vanish for $2k-1>\text{corank}(\F)$.

We point out that these classes were known to people working in foliation 
theory (see e.g.~\cite{Goldman:article:1}).

\paragraph{Poisson Manifolds}
Let $(M,\Pi)$ be a Poisson manifold, so $A=T^*M$. A basic connection
on $A$ is a basic contravariant connection $\nabla$ in the sense of
\cite{Fernandes:article:1}. and we can take for $\nabla^1$ on $E=T^*M\ominus
T^*M$ the connection $\nabla^1=\nabla\oplus\nabla$. Also, we let
$\bar{\nabla}$ be some riemannian connection on $T^*M$ and set
$\nabla^0=\bar{\nabla}\oplus\bar{\nabla}$. It is clear that
\[ \lambda(\nabla^1,\nabla^0)(P)=2\lambda(\nabla,\bar{\nabla})(P),\]
and it follows that the characteristic classes we have defined for
$A=T^*M$ are equal to twice the characteristic classes we have defined
in \cite{Fernandes:article:1} for the case of a Poisson manifold:
\[ m_k(T^*M)=2m_k(M).\]

A special case where computations can be made explicitly is when
$M=\gg^*$ with the Lie-Poisson structure. We have shown in
\cite{Fernandes:article:1} that these classes are represented by
Lie algebra cocycles given by the general formula:
\begin{align}
\label{characteristic:dual:Lie:algebra}
m_k(\gg^*)(v_1,\dots,v_{2k-1})&=\nonumber\\
\frac{1}{(2\pi)^k}\sum_{\sigma\in S_{2k-1}}&
K_k(v_{\sigma(1)},[v_{\sigma(2)},v_{\sigma(3)}],\dots,
[v_{\sigma(2k-2)},v_{\sigma(2k-1)}])
\end{align}
where:
\[ K_j(v_1,\dots,v_j)\equiv\tr(\ad v_1\cdots\ad v_j).\]

\paragraph{Transformation Lie algebroids}
Consider an infinitesimal action $\rho:\gg\to \X^1(M)$ of a Lie
algebra $\gg$ on a manifold $M$, so $A=M\times\gg$. Sections of $A$
can be identified with $\gg$-valued functions on $M$, so if $v\in \gg$
we identify it with a constant section.

There is a canonical choice of basic $A$-connection on $A$, namely the
unique $A$-connection which for constant sections satisfies:
\[ \nabla_v w=[v,w].\]
This connection is compatible with the Lie algebroid since
$\#\nabla_v w=\check{\nabla}_v\#w$ where $\check{\nabla}$ is the
$A$-connection on $TM$ which for any constant section $v\in\gg$ and
vector field $X\in\X^1(M)$ satisfies
\[ \nabla_{v}X=\Lie_{\rho(v)}X.\]
Hence we take the connection $\nabla^1$ on $E=A\ominus T^*M$ given
by
\[\nabla^1_v (w+\omega)=[v,w]-\Lie_{\rho(v)}\omega.\]

Now pick some riemannian metric on $M$ and consider the flat metric on
$\gg$. This gives a riemannian connection $\nabla^0$ on $E$ and we can
then compute the secondary characteristic classes. Formula
(\ref{eq:invariants:2}) shows that the classes
$\lambda^{1,0}(P)$, in general, will depend on
the curvature $R^t$ of $t\nabla^1+(1-t)\nabla^0$ in some intricate
way. However, for the first characteristic class where
$P=\frac{1}{2\pi}\tr$, there is no dependence on the curvature, and we
find explicitly:
\[m_1(A)(x)=\frac{1}{2\pi}\tr\ad x-\div(\rho(x))\]
(again, we view $x\in\gg$ has a constant section of $A$), where
$\div$ is the divergence operator on vector fields defined by the
metric on $M$.

Assume further that $M=V$ is a vector space and $\rho$ is a linear
action (a Lie algebra representation). In this case we can also choose
on $V$ a flat metric and it is possible to compute all characteristic
classes. Let $\tilde{\rho}$ denote the direct sum of representations
$\ad\oplus \rho^*$ on $\gg\times V^*$, where
$\rho^*:\gg\to\mathfrak{gl}(V^*)$ is the dual representation to $\rho
:\gg\to\mathfrak{gl}(V)$. Then the general formula for the
characteristic classes $m_k$ is:
\begin{multline}
m_k(A)(x_1,\dots,x_{2k-1})=\\ C_k \sum_{\sigma\in S_{2k-1}}
P_k(\tilde{\rho}(v_{\sigma(1)}),\tilde{\rho}([v_{\sigma(2)},v_{\sigma(3)}]),
\dots,
\tilde{\rho}([v_{\sigma(2k-2)},v_{\sigma(2k-1)}]))
\end{multline}
where $P_k$ are the elementary symmetric polynomials and $C_k$ is
some numerical factor.

For example, if we let 
\[\rho=\ad:\gg\to\mathfrak{gl}(\gg^*)\]
be the
coadjoint action of $\gg$ we see that $m_k(A)$ is twice the Lie
algebra cohomology classes given by
(\ref{characteristic:dual:Lie:algebra}), i.e., they coincide with
the characteristic classes for $A=T^*M$, where we take $M=\gg^*$
with the Lie-Poisson bracket.

\begin{acknowledgment}
I would like to thank Alan Weinstein for suggesting Lie algebroids
as the proper setting for developing a full theory of connections
and characteristic classes, as well as for several illuminating
discussions. His joined paper \cite{Evens:article:1} with Sam Evens 
and Jiang-Hua Lu was a source of inspiration for the present 
work. I also would like to thank Marius Crainic and Victor L.~Ginzburg for 
additional discussions, remarks and suggestions, and the anonymous 
referee for enumerous observations that have improved the manuscript.
\end{acknowledgment}

\end{article}
\end{document}